\documentclass[preprint,11pt]{elsarticle}

\usepackage{amsmath,amssymb,latexsym}
\usepackage{booktabs}
\usepackage{url}
\usepackage{color}
\journal{Linear Algebra and Its Applications}

\usepackage[mathscr]{euscript}

\newtheorem{theorem}{Theorem}[section]

\newtheorem{lemma}[theorem]{Lemma}

\newtheorem{corollary}[theorem]{Corollary} 

\newtheorem{definition}[theorem]{Definition}
\newtheorem{remark}[theorem]{Remark}

\newproof{proof}{Proof}

\usepackage{stmaryrd}
\usepackage{tikz}
\usepackage{booktabs}

\newcommand{\RS}[2]{\textrm{RS}(#1\kern-2pt\mid\kern-2pt#2)}
\newcommand{\RSM}[2]{(#1\kern-2pt\mid\kern-2pt#2)}
\newcommand{\smatrix}[4]{%
  \left(\begin{smallmatrix}#1&#2\\#3&#4\end{smallmatrix}\right)}

\newcommand{\F}{{\mathbb F}}
\newcommand{\GL}{{\mathrm{GL}}}
\newcommand{\SL}{{\mathrm{SL}}}
\newcommand{\OmegaX}{\Omega\mathbf{X}}
\newcommand{\GX}{\mathrm{G}\mathbf{X}}
\newcommand{\coloneq}{:=}
\newcommand{\SX}{{\mathrm{S}\mathbf{X}}}
\newcommand{\SY}{{\mathrm{S}\mathbf{Y}}}
\newcommand{\SO}{{\mathrm{SO}}}

\newcommand{\lhdeq}{\trianglelefteqslant}    
\newcommand{\Z}{{\mathbb Z}}
\newcommand{\W}{X_1\kern-1pt\times\kern-1pt X_2\kern-1pt\times\kern-1pt X_1\kern-1pt\times\kern-1pt X_2}
\newcommand{\WW}{X_2\kern-1pt\times\kern-1pt X_1\kern-1pt\times\kern-1pt X_2\kern-1pt\times\kern-1pt X_1}
\newcommand{\cC}{{\mathscr{C}}}
\newcommand{\cR}{R_{e_2\times e_1}}  
\newcommand{\cY}{{\mathcal{Y}}}
\newcommand{\cE}{{\mathscr{E}}}

\numberwithin{equation}{section}

\begin{document}
\title{Bipartite $q$-Kneser graphs and\\ two-generated irreducible linear groups}

\author{S.P. Glasby}
\ead{Stephen.Glasby@uwa.edu.au}
\affiliation[1]{
	organization={Centre for the Mathematics of Symmetry and Computation,
        The University of Western Australia},
	addressline={35 Stirling Highway},
	postcode={6009},
	city={Perth},
	country={Australia}
}
\author{Alice C. Niemeyer}
\ead{Alice.Niemeyer@art.rwth-aachen.de}
\affiliation[2]{
	organization={Chair for Algebra and Representation Theory,
	 RWTH Aachen University},
	addressline={Pontdriesch 10-16},
	postcode={52062},
	city={Aachen},
	country={Germany}
        }
\author{Cheryl E. Praeger}

\ead{Cheryl.Praeger@uwa.edu.au}  
\affiliation[3]{
	organization={Centre for the Mathematics of Symmetry and Computation,
        The University of Western Australia},
	addressline={35 Stirling Highway},
	postcode={6009},
	city={Perth},
	country={Australia}
}


\begin{keyword}
$q$-Kneser graph\sep irreducible linear group\sep $k$-walk\sep
      closed $k$-arc\sep stingray element
\MSC[2020] Primary: 05C50, 20-08; Secondary: 05C35.
\end{keyword}

\begin{abstract}
  Let $V\coloneq(\F_q)^d$ be a $d$-dimensional vector space over the
  field $\F_q$ of order~$q$. Fix positive integers $e_1,e_2$ satisfying
  $e_1+e_2=d$.
  Motivated by analysing a fundamental algorithm in computational group
  theory for recognising classical groups, we consider a certain
  quantity $P(e_1,e_2)$ which arises in both graph theory and
  group representation theory: $P(e_1,e_2)$ is the proportion of
  $3$-walks in the `bipartite $q$-Kneser graph' $\Gamma_{e_1,e_2}$ that are
  closed $3$-arcs. We prove that, for a  group $G$ satisfying $\SL_d(q)\lhdeq G\leqslant \GL_d(q)$,
  the proportion of certain element-pairs  in $G$ called `$(e_1,e_2)$-stingray duos' which generate an irreducible subgroup is also equal to $P(e_1,e_2)$.
  We give an exact formula for $P(e_1,e_2)$, and prove that
  \[1-q^{-1}-q^{-2}< P(e_1,e_2)< 1-q^{-1}-q^{-2}+2q^{-3}-2q^{-5}\]
  for $2\leqslant e_2\leqslant e_1$ and $q\geqslant2$.
 These bounds have implications for the complexity analysis of the state-of-the-art
 algorithms to    recognise classical groups, which we discuss in the final section.
\end{abstract}

\maketitle


\section{Introduction}\label{sec:intro}

  Let $V\coloneq(\F_q)^d$ be a $d$-dimensional vector space over the
  field $\F_q$ of order~$q$. Given positive integers $e_1,e_2$ satisfying
  $e_1+e_2=d$, let $X_i$ denote the set of $e_i$-dimensional subspaces of~$V$
  for $i\in\{1,2\}$. The \emph{bipartite $q$-Kneser graph}
  $\Gamma=\Gamma_{e_1,e_2}$ has vertex set
  $V\Gamma=X_1\dot{\cup} X_2$ (disjoint union) and edge set
  $E\Gamma=\{\{S_1,S_2\} \mid S_1\in X_1, S_2\in X_2,   S_1\cap S_2=0\}$.
  We emphasise that the subsets $X_1, X_2$ are regarded as disjoint vertex-subsets, even in the case $e_1=e_2$.
  A $4$-tuple
  $(S_0,S_1,S_2,S_3)\in \WW$ is called a
  \emph{$3$-walk} if $\{S_0,S_1\}$, $\{S_1,S_2\}$, $\{S_2,S_3\}$ are edges
  of~$\Gamma$; a $3$-walk is called \emph{closed}
  if $\{S_3,S_0\}$ is also an edge, and it is called
  a \emph{$3$-arc} if $S_0\ne S_2$ and $S_1\ne S_3$. Define $P(e_1,e_2)$ to be
  the proportion of $3$-walks of $\Gamma_{e_1,e_2}$ which are closed $3$-arcs.
  Our first result determines an exact formula for $P(e_1,e_2)$ in terms of the 
 following function:
  \begin{equation}\label{e:om}
    \omega(e)\coloneq\prod_{i=1}^e(1-q^{-i}),
    \qquad \mbox{for $e\in\Z$ with $e\geqslant1$, and $\omega(0)=1$.}
    \end{equation}

  \begin{theorem}\label{TT:prop}
    Let $d=e_1+e_2$ with $1\leqslant e_2\leqslant e_1$.
    Let $q>1$ be a prime power, and let
    $P(e_1,e_2)$ denote the proportion of $3$-walks
    of $\Gamma_{e_1,e_2}$ which are closed $3$-arcs. Then   
  \begin{align}\label{e:Pe12}
  P(e_1,e_2)&=-(1-q^{-e_1e_2})q^{-e_1e_2}+\sum_{\ell=0}^{e_2-1}
  \frac{\omega(e_1)\omega(e_2)q^{-(e_1-e_2+\ell)\ell}}{\omega(e_1-e_2+\ell)\omega(\ell)}\\
  &=1-O(1/q).\notag
  \end{align}
  \end{theorem}

  Although the exact formula~\eqref{e:Pe12} for the proportion $P(e_1,e_2)$
  is intricate, it allows us to prove that $P(e_1,e_2)$ is
  close to $1-q^{-1}-q^{-2}$. This requires very delicate calculations,
  especially for `small' $q$. The lower bound given was our primary
  objective to help understand the complexity of a probabilistic
  generation algorithm, especially for `small'~$q$.

  \begin{theorem}\label{TT:bounds}
  If $e_2=1$, then 
  \begin{align*}
  &1-q^{-1}-q^{-2}< P(e_1,1)<1-q^{-1} \quad\textup{for $e_1\geqslant3$, and}\\
  &1-q^{-1}-q^{-2}< P(e_1,e_2)<1-q^{-1}-q^{-2}+2q^{-3}-2q^{-5}
  \quad\textup{for $2\leqslant e_2\leqslant e_1$.}
  \end{align*}
  \end{theorem}

  The proportion $P(e_1,e_2)$ also arises when considering
  certain 2-generated irreducible subgroups of $\GL(V)$.
  To describe this key connection, we need additional terminology.

  An element $g\ne1$ of the general linear group $\GL(V)=\GL_d(q)$ is called an
  \emph{$e$-stingray} element if $g$ acts irreducibly on the image
  $U=V(g-1)={\rm im}(g-1)$ of $g-1$, and $\dim(U)=e$. A pair
  $(g_1,g_2)$ of elements in $\GL(V)$ is called an
  \emph{$(e_1,e_2)$-stingray duo} if $g_i$ is an $e_i$-stingray
  element for $i\in\{1,2\}$, and $U_1\cap U_2=\{0\}$ where
  $U_i={\rm im}(g_i-1)$. 
  A pair $(g_1,g_2)\in G\times G$ is called
  \emph{irreducible} if the 2-generated subgroup $\langle g_1,g_2\rangle$ of
  $\GL(V)$ acts irreducibly, that is, the only subspaces of $V$ invariant
  under $\langle g_1,g_2\rangle$ are $V$ and $\{0\}$. 
  Clearly, for an $(e_1,e_2)$-stingray duo  $(g_1,g_2)$, the subgroup
  $\langle g_1,g_2\rangle$ is reducible whenever $e_1+e_2<d$
  as it fixes the proper subspace $U_1+U_2$. The graph-theoretic
  invariant $P(e_1,e_2)$ described in Theorem~\ref{TT:prop} has a
  representation-theoretic interpretation as described below.

  \begin{theorem}\label{TT:stingray}
  If $d=e_1+e_2$ with $1\leqslant e_2\leqslant e_1$, and if $\SL_d(q)\trianglelefteqslant G\leqslant \GL_d(q)$, 
  then with $P(e_1,e_2)$ as in \eqref{e:Pe12},
\[
P(e_1,e_2)=\frac{\text{Number of irreducible $(e_1,e_2)$-stingray duos in $G\times G$}}{\text{Number of $(e_1,e_2)$-stingray duos in $G\times G$}}.
\]  
\end{theorem}

  An $(e_1,e_2)$-stingray duo $(g_1,g_2)$ is called \emph{generating} if
  $\langle g_1,g_2\rangle$ contains the special linear group $\SL_d(q)$.
  For algorithmic purposes we are interested in the proportion of
  $(e_1,e_2)$-stingray duos that are generating. This will be a smaller
  proportion than the proportion $P(e_1,e_2)$ of stingray duos that are
  irreducible.
  We believe that this smaller proportion is also $1-O(q^{-1})$,
  but a proof of this fact requires a very careful analysis of when
  $\langle g_1,g_2\rangle\cap\SL_d(q)$ is a proper subgroup of $\SL_d(q)$.
  When non-generation occurs, it is almost always because
  $\langle g_1,g_2\rangle$ is a reducible subgroup. By the
  previous two theorems, the proportion of
  $(e_1,e_2)$-stingray duos in $G\times G$ that are reducible is $1-P(e_1,e_2)< q^{-1}+q^{-2}$.
    By contrast, we believe that the proportion of
  $(e_1,e_2)$-stingray duos that are non-generating and irreducible is
  substantially smaller:
  at most $O(q^{-ce_1e_2})$ for some constant $c$. 
   This is known to be the case when $e_1=e_2=d/2$, see \cite[Theorems 5 and 6]{PSY} and  Section~\ref{S:comp}, which is the case underpinning  the algorithm presented in \cite{DLLOB}.
    Theorem~\ref{TT:bounds} and Corollary~\ref{C:DLLOB} lead to an
    improved complexity analysis for the algorithm
    in~\cite{DLLOB}, for further details see the forthcoming paper~\cite{GNP2023}.
    

The estimates obtained in this paper will be applied in \cite{GNP2023} to prove a key theorem underpinning the complexity analysis of a new generation of recognition algorithms for classical groups.
These new algorithms are described in Rademacher's PhD thesis \cite{DR} where, in addition, information is given on implementation details and 
comparative timings for the new algorithm against the current-best algorithm in \cite{DLLOB} for various values of the dimension and field size. The latter indicate considerably  improved running times for the new algorithm, especially as the dimension increases.

  The bipartite $q$-Kneser graph $\Gamma_{e_1,e_2}$ was also used to solve
  another problem that arises from computational group theory. Suppose
  that $V$ is a $(e_1+e_2)$-dimensional classical space (symplectic or
  orthogonal) over a finite field $\F_q$, and let $Y_i$ be the set of
  \emph{non-degenerate} $e_i$-subspaces of~$V$ (of a particular type
  in the orthogonal case).  A pair $(U_1,U_2)\in Y_1\times Y_2$ is
  called \emph{spanning} if $V=U_1+U_2$, that is to say, if $\{U_1,
  U_2\}$ is an edge of $\Gamma_{e_1,e_2}$. Thus the proportion of
  pairs $(U_1,U_2)\in Y_1\times Y_2$ that are spanning equals the
  proportion of pairs $(U_1,U_2)\in Y_1\times Y_2$ that are edges of
  the induced subgraph $[Y_1\,\dot{\cup}\, Y_2]$ of $\Gamma_{e_1,e_2}$ with
  vertex set $Y_1\,\dot{\cup}\, Y_2$. It was shown in~\cite[Theorem 1.1]{GNP2022a} that this proportion is $1-O(q^{-1})$, and a stronger estimate was obtained in~\cite[Theorem~1.1]{GIM}: namely the proportion is at most $1-\frac{3}{2q}$ unless
  $(e_1,e_2,q)=(1,1,2)$. (The bound also applies for a unitary space
  $(\F_{q^2})^d$ if we replace $q$ with~$q^2$.) The proof in \cite{GIM} uses the Bipartite Expander
  Mixing Lemma, and knowledge of the eigenvalues of the induced
  subgraph $[Y_1\,\dot{\cup}\, Y_2]$ which are determined via a deep
  geometric algorithm~\cite{DMM} dating back to
  Brouwer~\cite{Brouwer2010}, and the representation theory of the
  symmetric group. By contrast, the present paper uses elementary
  combinatorics, linear algebra and group theory to prove the
  above~theorems.


  The research in this paper sheds light on a very difficult part
  of a larger problem of recognising
  classical groups computationally.
  Much has been written on the computational complexity of classical recognition algorithms 
  including~\cite{DLLOB,GNP2022a,GNP2022b,LOB,KS,PSY}, and with the benefit
  of hindsight, the `right' proportions were not always studied.
  An overview of this larger picture, and some of the nuances, is discussed
  more fully in Section~\ref{S:comp}. In particular,
  by noting that some of these algorithms work with stingray duos instead of `stingray pairs' we obtain probability estimates which allow the hypothesis $q>4$   to be removed from~\cite[Theorem~1.2]{DLLOB}, and  the hypothesis $q>2$ to be removed from~\cite[Theorem~2]{PSY}. Application of our  Theorem~\ref{TT:bounds} also improves the analogous upper bound for duos in the linear case in \cite[Theorem~5]{PSY}.

  Although a major motivation for this paper is analysing a classical recognition algorithm, including the case of special linear groups, we stress that this paper is \emph{motivated by the linear case}. The hard problem for the analysis is showing that the subgroup $\langle g_1,g_2\rangle$, generated by a stingray duo $(g_1,g_2)$, lies inside a proper subgroup of the quasisimple classical group, with low probability. It is shown in~\cite{GNP2023} that in the symplectic, unitary, and orthogonal cases the
  probability is at most $c_1q^{-c_2d}$ for
  positive constants $c_1,c_2$. By contrast,
  in the linear case, Theorems~\ref{TT:bounds}  and~\ref{TT:stingray} of this paper imply that even the probability that $\langle g_1,g_2\rangle$ lies in a reducible subgroup of $\GL_d(q)$ is much higher, namely just less than $q^{-1}+q^{-2}$. 
  Obtaining this highly accurate estimate has proved critical for giving complexity estimates for the recognition algorithm for $\SL_d(q)$ that are effective even for small values of $q$ (which previous results were unable to handle).


    In Section~\ref{S:Kneser}, we give a more general definition of a
  bipartite $q$-Kneser graph where $e_1+e_2\leqslant d$, but because of our applications, we focus on the case where $e_1+e_2=d$. We count the number of
  $3$-walks and $3$-arcs in $\Gamma_{e_1,e_2}$, see Remark~\ref{R:WW} and Lemma~\ref{L:number}.
  In Section~\ref{S:sting}, we analyse stingray pairs and duos. Lemma~\ref{L1}
  gives criteria for a stingray duo $(g_1,g_2)$ to be irreducible in terms of subspaces $U_1=\mathrm{im}(g_1-1)$ and $U_2=\mathrm{im}(g_2-1)$, and these are shown in Lemma~\ref{L:irredclosed} to hold 
  precisely when we have a closed 3-arc in $\Gamma_{e_1,e_2}$.  Section~\ref{S:3-walks} counts the number of
  closed 3-walks and closed 3-arcs of $\Gamma_{e_1,e_2}$.
  These counts are used in Section~\ref{S:3-arcs} to compute the proportion
  $P(e_1,e_2)$ of 3-walks that are closed 3-arcs. The  $q$-identity in Corollary~\ref{C} allows us to highlight the dominant terms of
  $P(e_1,e_2)$. This is used in Section~\ref{S:Order} to prove
  the upper and lower bounds for $P(e_1,e_2)$ in Theorem~\ref{TT:bounds}. Finally, the computational context for this work is
  described in Section~\ref{S:comp}.

\section{The {$q$}-Kneser graph}\label{S:Kneser}

\subsection{Walks and arcs in graphs}
A \emph{$k$-walk} of a graph $\Gamma$ is a sequence $(v_0,v_1,\dots,v_k)$
of $k+1$ vertices of $\Gamma$ such that $\{v_{i-1},v_i\}$ is an edge
for $1\leqslant i\leqslant k$. We call a $k$-walk  $(v_0,v_1,\dots,v_k)$  
of $\Gamma$ a \emph{$k$-arc} if $v_{i-1}\ne v_{i+1}$ for $1\leqslant i<k$,
see~\cite[p.\;130]{Biggs}. For a $k$-arc, imagine `walking' from $v_0$ to $v_k$
via edges of $\Gamma$, where we cannot walk from $v_i$ directly back to
$v_{i-1}$, but we can walk to any other adjacent vertex.
A 1-arc $(v_0,v_1)$ is commonly called an \emph{arc}, and is viewed
as the edge $\{v_0,v_1\}$ directed from $v_0$ to~$v_1$.  A $k$-walk or a $k$-arc $(v_0,v_1,\dots,v_k)$
is called \emph{closed}\footnote{Our definition of `closed $k$-walk' differs from~\cite{Biggs}. The $(k+2)$-tuple $(S_0,\dots,S_k,S_0)$ is a closed \emph{$(k+1)$-walk} according to Biggs~\cite[p.\,12]{Biggs}. For us, the $(k+1)$-tuple $(S_0,\dots,S_k)$ a closed \emph{$k$-walk} when $\{S_k,S_0\}$ is~an~edge.} if $\{v_k,v_0\}$ is an edge of $\Gamma$. If $\Gamma$ has
$n$ vertices and its adjacency matrix has eigenvalues
$\lambda_1,\dots, \lambda_n$ (counting multiplicities), then it follows
from~\cite[Additional Results~2h]{Biggs} that the number of 
closed $k$-walks of $\Gamma$ equals $\sum_{i=1}^n\lambda_i^{k+1}$.
We define the bipartite $q$-Kneser graph below, and show that the number of closed
3-walks can be counted \emph{without} using eigenvalues.

If $\Gamma$ is
a bipartite graph, then a 3-arc has four
\emph{distinct} vertices, as the only edges join the `left' vertices to
the `right' vertices of $\Gamma$.

\subsection{Bipartite $q$-Kneser graphs}

Let $V\coloneq(\F_q)^d$
be a $d$-dimensional vector space over the field $\F_q$ of order~$q$.
Given positive integers $e_1,e_2$ satisfying $e_1+e_2\leqslant d$, let
\begin{equation}\label{e:x12}
\mbox{$X_i$ denote the set of $e_i$-subspaces of $(\F_q)^d$ for $i\in\{1,2\}$.}
\end{equation} 
We follow~\cite{GIM} and define the \emph{bipartite $q$-Kneser graph}
$\Gamma\coloneq\Gamma_{d,e_1,e_2}$. The vertex set $V\Gamma$ is the disjoint union
$X_1\dot{\cup} X_2$, and the edge set $E\Gamma$ comprises all 2-subsets
$\{S_1,S_2\}$ with $S_1\in X_1$, $S_2\in X_2$ and $S_1\cap S_2=\{0\}$.
If $e_1+e_2>d$, then $\Gamma_{d,e_1,e_2}$ has no edges, and if $d=e_1+e_2$,
then we write $\Gamma_{e_1, e_2}$ instead of $\Gamma_{d,e_1, e_2}$.

When $e_1=e_2$ there is a \emph{non-bipartite $q$-Kneser graph}
$\widetilde{\Gamma}$ with vertex set $X_1=X_2$ and $\{S_1,S_2\}$ is an
edge precisely when $S_1\cap S_2=\{0\}$.
The bipartite graph $\Gamma$ is the standard bipartite double-cover of the
non-bipartite graph $\widetilde{\Gamma}$,
so that $\lambda$ is an eigenvalue of $\widetilde{\Gamma}$ if and only
if $\pm\sqrt{\lambda}$ are eigenvalues of $\Gamma$
as explained in~\cite{GIM}. We henceforth consider $\Gamma$ and
not~$\widetilde{\Gamma}$. We warn the reader that we sometimes
refer to $\Gamma$ simply as the \emph{$q$-Kneser graph},
omitting the important adjective~`bipartite'.

Because of a link which we uncover in Section~\ref{S:sting} between
irreducible $2$-generated linear groups and closed $3$-arcs in $\Gamma_{e_1, e_2}$
(see Lemma~\ref{L:irredclosed}) we shall assume after Section~\ref {S:3-walks}
that $e_1+e_2=d$ holds.
We refer to the quantity $\omega(e,q)=|\GL_e(q)|/q^{e^2}$
simply as $\omega(e)$, suppressing $q$, as per the definition
given in \eqref{e:om}.
We count the number of 3-walks  of $\Gamma_{e_1+e_2,e_1,e_2}$, and
the number of 3-arcs, in Lemma~\ref{L:number}. The
different types of $3$-walks, with possible repeated vertices, are
illustrated in Figure~\ref{F1}; only the leftmost represents a $3$-arc. 

\begin{figure}[!ht]
\begin{center}
\begin{tikzpicture}[scale=0.6] 
\node [left] at (0.7,4.1) {$X_1$};
\node [right] at (1.5,4.1) {$X_2$};
\node [left] at (-0.2,3) {$U_1$};
\node [left] at (-0.2,1) {$F_2$};
\node [right] at (2.2,1) {$F_1$};
\node [right] at (2.2,3) {$U_2$};
\draw (0,1)--(2,3)--(0,3)--(2,1);
\draw [thick,dashed] (0,1)--(2,1);
\draw [fill] (0,1) circle [radius=0.08];
\draw [fill] (2,1) circle [radius=0.08];
\draw [fill] (0,3) circle [radius=0.08];
\draw [fill] (2,3) circle [radius=0.08];
\draw (0,2) circle [x radius=0.35, y radius=1.8];
\draw (2,2) circle [x radius=0.35, y radius=1.8];
\end{tikzpicture}
\hskip3mm 
\begin{tikzpicture}[scale=0.6] 
\node [left] at (0.7,4.1) {$X_1$};
\node [right] at (1.5,4.1) {$X_2$};
\node [left] at (-0.2,2) {$U_1$};
\node [right] at (2.2,1) {$F_1$};
\node [right] at (2.2,3) {$U_2$};
\draw (0,2)--(2,3)--(0,2)--(2,1);
\draw [fill] (2,1) circle [radius=0.08];
\draw [fill] (0,2) circle [radius=0.08];
\draw [fill] (2,3) circle [radius=0.08];
\draw (0,2) circle [x radius=0.35, y radius=1.8];
\draw (2,2) circle [x radius=0.35, y radius=1.8];
\end{tikzpicture}
\hskip3mm 
\begin{tikzpicture}[scale=0.6] 
\node [left] at (0.7,4.1) {$X_1$};
\node [right] at (1.5,4.1) {$X_2$};
\node [left] at (-0.2,3) {$U_1$};
\node [left] at (-0.2,1) {$F_2$};
\node [right] at (2.2,2) {$U_2$};
\draw (0,1)--(2,2)--(0,3)--(2,2);
\draw [fill] (0,1) circle [radius=0.08];
\draw [fill] (2,2) circle [radius=0.08];
\draw [fill] (0,3) circle [radius=0.08];
\draw (0,2) circle [x radius=0.35, y radius=1.8];
\draw (2,2) circle [x radius=0.35, y radius=1.8];
\end{tikzpicture}
\hskip3mm 
\begin{tikzpicture}[scale=0.6] 
\node [left] at (0.7,4.1) {$X_1$};
\node [right] at (1.5,4.1) {$X_2$};
\node [left] at (-0.2,2) {$U_1$};
\node [right] at (2.2,2) {$U_2$};
\draw (0,2)--(2,2);
\draw [fill] (0,2) circle [radius=0.08];
\draw [fill] (2,2) circle [radius=0.08];
\draw (0,2) circle [x radius=0.35, y radius=1.8];
\draw (2,2) circle [x radius=0.35, y radius=1.8];
\end{tikzpicture}
\end{center}
\vskip-5mm
\caption{3-walks $(F_1,U_1,U_2,F_2)$ with possible repeated vertices; in the second, third and fourth diagrams we have $U_1=F_2$, $U_2=F_1$ and $(U_1,U_2)=(F_2,F_1)$, respectively}
\label{F1}
\end{figure}

\begin{remark}\label{R:WW}
  As $\Gamma_{d,e_1,e_2}$ is a bipartite graph with vertex set $X_1\dot{\cup} X_2$,
  there are two types of 3-walks: those in
  $W_1\coloneq \W$ and those in
  $W_2\coloneq \WW$. The reversal map
  $(S_0,S_1,S_2,S_3)\mapsto(S_3,S_2,S_1,S_0)$ is a bijection $W_1\to W_2$
  which preserves 3-walks, 3-arcs, closed 3-walks, and closed 3-arcs.
  Hence the \emph{proportion} of 3-walks (resp. 3-arcs) in $\Gamma_{d,e_1,e_2}$ that are
  closed equals the proportion of 3-walks (resp. 3-arcs) in $W_2$ that are
  closed. This explains why we restrict to $W_2$ in Lemma~\ref{L:number} and Theorems~\ref{T:walk},~\ref{T:3arcs}.
\end{remark}

\begin{lemma}\label{L:number}
  If $d=e_1+e_2$, then the number of $3$-walks in $\WW$ is $q^{4e_1e_2}\xi$,
  and the number of $3$-arcs is $q^{4e_1e_2}\xi(1-q^{-e_1e_2})^2$,
  where $\xi=\frac{\omega(e_1+e_2)}{\omega(e_1)\omega(e_2)}$ with $\omega(e)$ as in \eqref{e:om}.
\end{lemma}

\begin{proof}
  We first count the number of $3$-walks
  $(F_1,U_1,U_2,F_2)\in \WW$.
  The number of choices of $F_1$ is
  \begin{equation}\label{E:w}
    |X_2|=|\textup{$\{\,e_2$-subspaces of $(\F_q)^{e_1+e_2}\,\}$}| = q^{e_1e_2}\xi
    \quad\textup{where}\quad
    \xi\coloneq\frac{\omega(e_1+e_2)}{\omega(e_1)\omega(e_2)}.
  \end{equation}
  There are $q^{e_1e_2}$ choices for each of the complements $U_1$ of $F_1$, $U_2$
  of $U_1$ and $F_2$ of $U_2$, so the number of 3-walks of $\Gamma$ is
  $q^{4e_1e_2}\xi$ with $\xi$ as in the statement. Similarly,  the number
  of 3-arcs of $\Gamma$ is
  $q^{e_1e_2}\xi q^{e_1e_2}(q^{e_1e_2}-1)^2=q^{4e_1e_2}\xi(1-q^{-e_1e_2})^2$
  as $U_2\ne F_1$ and $F_2\ne U_1$.
\end{proof}

\begin{remark}
  For $d=e_1+e_2$ and $e_2\leqslant e_1$ it is proved in~\cite{GIM} that
  $\Gamma_{e_1,e_2}$ has $2(e_2+1)$ distinct eigenvalues and their values are
  $\pm\mu_0,\dots,\pm\mu_{e_2}$ where $\mu_j=q^{e_1e_2-j(e_1+e_2-j)/2}$
  for $0\leqslant j\leqslant e_2$. Let $M_j$ denote the \emph{multiplicity} of
  the eigenvalue $\mu_j$. This is also the multiplicity of $-\mu_j$.
  By~\cite[Additional Results~2h]{Biggs} the number of closed $k$-walks
  is $2\sum_{j=0}^{e_2} M_j \mu_j^{k+1}$ and hence the number of closed 3-walks
  in $\Gamma_{e_1,e_2}$ is $2\sum_{j=0}^{e_2} M_j q^{4e_1e_2-2j(e_1+e_2-j)}$.
  Computing formulas for the $\mu_j$ is not elementary. Indeed, the proof
  in~\cite{GIM} uses the representation theory of the symmetric group, and
  a sophisticated geometric algorithm due to~\cite{Brouwer2010} and~\cite{DMM}.
  We shall compute both the number of closed 3-walks in $\Gamma_{e_1,e_2}$,
  and the number of closed 3-arcs,  by using linear algebra and
  group theory only. This bypasses the need to calculate the multiplicities
  $M_j$, $0\leqslant j\leqslant e_2$, which were not described
  in~\cite{GIM} for $j\ne0$. Our proof
  involves elementary arguments only.\qed
\end{remark}

\section{Stingray duos and the $q$-Kneser graph}\label{S:sting}

We view $V\coloneq(\F_q)^d$ as the natural module for the general linear group
$G_d\coloneq\GL_d(q)$. 
Let $e_1+e_2\leqslant d$ with each $e_i\geqslant 1$.
Fix $g_1,g_2\in\GL_d(q)$ and set $U_i\coloneq\textup{im}(g_i-1)$,
$F_i\coloneq\textup{ker}(g_i-1)$ and $e_i=\dim(U_i)$ for $i\in\{1,2\}$
where $\textup{ker}(x)$ denotes the kernel of a linear transformation $x$, and
as in Section~\ref{sec:intro}, $Vx=\textup{im}(x)$ denotes its image.

\begin{definition}\label{d:sting}
  An element $g\in \GL(V)$ is called a \emph{stingray element} if $g$ acts
  irreducibly and non-trivially on $U\coloneq\textup{im}(g-1)$.
A stingray element $g$ is called an \emph{$e$-stingray element}
if $\dim(U)=e$.
(The first author coined this term as the matrix of $g$ looks like
\begin{tikzpicture}[scale=0.2] 
\draw [semithick] (0,1)--(0,2)--(1,2)--(1,1)--(0,1);
\draw [thick] (1,1)--(2,0);
\draw [fill, color=blue] (0.2,1.6) circle [radius=0.05];
\draw [fill, color=blue] (0.4,1.8) circle [radius=0.05];
\end{tikzpicture}
relative to a suitable basis, where the body has $\dim(U)$ rows, and
the tail has $\dim(V/U)$ rows.)
\end{definition}

\begin{definition}\label{d:duo}
  Given stingray elements $g_1,g_2\in\GL(V)$ set
  $U_i\coloneq\textup{im}(g_i-1)$ and $F_i\coloneq\textup{ker}(g_i-1)$ for
  $i\in\{1,2\}$.
  We call a pair $(g_1,g_2)$ a \emph{stingray pair}, or an
  $(e_1,e_2)$-\emph{stingray pair}, if $\dim(U_1)=e_1$ and $\dim(U_2)=e_2$.
  We call a stingray pair a \emph{stingray duo} if $U_1\cap U_2=\{0\}$.
\end{definition}

\begin{remark}\label{R}
  If the minimal polynomial of $g\in\GL(V)$ is a product $a(t)b(t)$ of coprime
  polynomials, then $V={\rm ker}(a(g))\oplus{\rm ker}(b(g))$ and
  ${\rm ker}(a(g))={\rm im}(b(g))$. If $g$ is an $e$-stingray element,  then its minimal polynomial has this form with $a(t)$ irreducible and $b(t)=t-1$, and we have
  $V=U\oplus F$ where $U={\rm ker}(a(g))$ and $F={\rm ker}(g-1)$.
  If $\dim(V)=e_1+e_2$, then an $(e_1,e_2)$-stingray duo has
  $V=U_1\oplus U_2=U_1\oplus F_1=U_2\oplus F_2$, and hence $\{U_1,U_2\}$,
  $\{U_1,F_1\}$ and $\{U_2,F_2\}$ are edges of $\Gamma_{e_1,e_2}$ and
  $(F_1,U_1,U_2,F_2)$ is a 3-walk of $\Gamma_{e_1,e_2}$.
  \qed 
\end{remark}

\begin{definition}\label{d:irredduo}
  A stingray duo $(g_1,g_2)$ is called an \emph{irreducible stingray duo}
  if $\langle g_1,g_2\rangle$ is an irreducible subgroup of $\GL(V)$,
  otherwise it is called a \emph{reducible stingray duo}.
\end{definition}

If $(g_1,g_2)$ is an $(e_1,e_2)$-stingray duo in $\GL(V)$, then
$e_1+e_2\leqslant d\coloneq\dim(V)$ since $U_1\oplus U_2\leqslant V$.
The following lemma is reminiscent
of~\cite[Lemma~3.7]{GNP2022b}. It characterises an irreducible stingray duo
$(g_1,g_2)$ in terms of the subspaces $U_i, F_i$. As $\langle g_1,g_2\rangle$ fixes $U_1+ U_2$, an
irreducible $(e_1,e_2)$-stingray duo must have $d=e_1+e_2$ and $V=U_1\oplus U_2$.

\begin{lemma}\label{L1}
  Suppose that $(g_1,g_2)$ is an $(e_1,e_2)$-stingray duo in $\GL_d(q)$.
  Then the $2$-generated subgroup $\langle g_1,g_2\rangle$ of $\GL_d(q)$
  acts irreducibly on $V=(\F_q)^d$ if and only if
  \begin{enumerate}[{\rm (a)}]
  \item $V=U_1\oplus U_2$ (so $d=e_1+e_2$),
  \item $F_1\cap F_2=\{0\}$, and
  \item $U_1\ne F_2$ and~$U_2\ne F_1$.
  \end{enumerate}
\end{lemma}

\begin{proof}
  Certainly $g_i$ preserves the decomposition $V=U_i\oplus F_i$, and
  fixes $F_i$ elementwise. Hence if $U_i\leqslant W\leqslant V$, then $g_i$ fixes
  $W$. Thus $U_1+U_2$ is invariant under $G\coloneq\langle g_1,g_2\rangle$,
  and also $G$ fixes $F_1\cap F_2$ elementwise. 
  Moreover, if $U_1= F_2$, then $G$ fixes $F_2$, and similarly
  if $U_2= F_1$, then $G$ fixes $F_1$. Hence if $G$ is irreducible,
  then the conditions (a)-(c) all hold.

  Conversely, suppose that (a)-(c) hold. We argue first that a $g_i$-invariant
  subspace $Z$ of $V$ satisfies $Z\leqslant F_i$ or $U_i\leqslant Z$. If some $z\in Z$
  satisfies $z\not\in F_i$, then $z(g_i-1)$ is a non-zero element of $U_i$.
  However, $g_i$ acts irreducibly on $U_i$ and so $U_i\leqslant Z$ as claimed.
  Suppose now that $W$ is a non-zero proper subspace invariant under
  $G\coloneq\langle g_1,g_2\rangle$. Then $W\leqslant F_1$ or $U_1\leqslant W$, and
  $W\leqslant F_2$ or $U_2\leqslant W$. It is not possible that $U_1\leqslant W$ and $U_2\leqslant W$,
  as then $V=U_1+ U_2\leqslant W$ by~(a). Also, it is not possible that $W\leqslant F_1$
  and $W\leqslant F_2$ as then $W\leqslant F_1\cap F_2=\{0\}$ by~(b).
  Hence either $W\leqslant F_1$ and $U_2\leqslant W$, or $W\leqslant F_2$ and $U_1\leqslant W$.
  In the former case, $U_2\leqslant F_1$ and $d=\dim(V)=e_1+e_2$ by (a), so
  that $e_2=\dim(U_2)\leqslant \dim(F_1)=d-e_1=e_2$ and hence $U_2=F_1$.
  In the latter case, similar reasoning shows that $U_1=F_2$. In either
  case, this contradicts~(c). Therefore $G$ preserves no proper non-zero
  subspace, and so $G$ acts irreducibly on $V$ as claimed.
\end{proof}

For  the subgroup $\langle g_1,g_2\rangle$ to be irreducible, we require $d=e_1+e_2$ by Lemma~\ref{L1}(a). Under this assumption we have multiple important links with the $q$-Kneser graph $\Gamma_{e_1,e_2}$.

\begin{lemma}\label{L:irredclosed}
  Let $d=e_1+e_2$. Let $g_1,g_2$ be stingray elements of $G_d=\GL_d(q)$.
  Write $e_i=\dim({\rm im}(g_i-1))$ and
  $\cC_i=\{x^{-1}g_ix\mid x \in G_d\}$ for $i\in\{1,2\}$.
\begin{enumerate}[{\rm (a)}]
\item For each $i$, the map $\phi_i: \cC_i\to E\Gamma_{e_1,e_2}\colon g\mapsto \{ {\rm im}(g-1), {\rm ker}(g-1)\}$ defines a surjection from $\cC_i$ to the
  set $E\Gamma_{e_1,e_2}$ of edges of $\Gamma_{e_1,e_2}$.
  \item Consider the map $\psi: \cC_1\times\cC_2 \to  \WW\colon (g_1',g_2')\mapsto (F'_1,U'_1,U'_2,F'_2)$, where  $U'_i={\rm im}(g_i'-1)$ and $F'_i={\rm ker}(g_i'-1)$, for $i=1,2$.
  \begin{enumerate}[{\rm (i)}]
    \item The pair $(g_1',g_2')$ is a stingray duo if and only if $\psi((g_1',g_2'))$ is a $3$-walk in $\Gamma_{e_1,e_2}$.
    \item Restricting $\psi$ to stingray duos yields a surjection onto the set of $3$-walks in~$\Gamma_{e_1,e_2}$.
    \item For an $(e_1,e_2)$-stingray duo $(g_1',g_2')$, the subgroup $\langle g_1',g_2'\rangle\leqslant \GL_d(q)$ acts irreducibly on $(\F_q)^d$ if and only if the image $\psi((g_1',g_2'))$ is a closed $3$-arc in $\Gamma_{e_1,e_2}$.
  \end{enumerate}	
\end{enumerate}	
\end{lemma}

\begin{proof}
(a) For $g\in\cC_i$ we see that $\{ {\rm im}(g-1), {\rm ker}(g-1)\}$
  forms a direct decomposition of the vector space $V=(\F_q)^d$, and hence is an edge of
  $\Gamma_{e_1,e_2}$ by Remark~\ref{R}. Each edge arises as an image
  of some element of $\cC_i$ under $\phi_i$ since $G_d=G_{e_1+e_2}$ acts
  transitively on the set of decompositions $V=U\oplus F$ with
  $\dim(U)=e_i$.

  (b) Note first that each conjugate of an $e$-stingray element is also an $e$-stingray
  element, and a conjugate in $G_d\times G_d$ of an
  $(e_1,e_2)$-stingray duo $(g_1,g_2)$ is also an $(e_1,e_2)$-stingray
  duo.  Fix $(g'_1,g'_2)\in\cC_1\times\cC_2$. Then, for each $i\in\{1,2\}$,
  we have $U'_i\cap  F'_i=\{0\}$ so $\{F'_1,U'_1\},\{U'_2,F'_2\}\in E\Gamma_{e_1,e_2}$. Thus
  $(F'_1,U'_1,U'_2,F'_2)$ is a 3-walk in $\Gamma_{e_1,e_2}$ if and only if
  $\{U'_1,U'_2\}$ is also an edge, that is, if and only if
  $U'_1\cap U'_2= \{0\}$. It follows from Definition~\ref{d:duo} that
  $(F'_1,U'_1,U'_2,F'_2)$ is a 3-walk if and only if $(g_1',g_2')$ is a
  stingray duo, proving part~(i). 
  
  Suppose that the $3$-walk
  $(F'_1,U'_1,U'_2,F'_2)$ is the image 
  of the stingray duo $(g_1',g_2')$ in $\cC_1\times \cC_2$ under $\psi$, and
  let $(F_1,U_1,U_2,F_2)$ be an arbitrary $3$-walk in~$\Gamma_{e_1,e_2}$. Since
  $G_{d}$ acts transitively on the decompositions $V=U'_i\oplus
  F'_i$ with $\dim(U'_i)=e_i$, for each $i=1,2$, there exist $x_1, x_2\in G_d$ such that $(F'_1,U'_1)^{x_1} = (F_1,U_1)$ and $(U'_2,F'_2)^{x_2}=(U_2,F_2)$. 
Then $\psi$ maps the pair $((g_1')^{x_1}, (g_2')^{x_2})\in \cC_1\times\cC_2$ to  $(F_1,U_1,U_2,F_2)$, and by part (i), $((g_1')^{x_1}, (g_2')^{x_2})$ is a stingray duo. This proves part~(ii).

  Finally, we prove part (iii). Let $(g_1',g_2')$ be a stingray duo, so
  $U'_i\cap F'_i=\{0\}$ and $\dim(F'_i)=d-e_i$ for each $i\in\{1,2\}$,
  and also $U'_1\cap U'_2= \{0\}$; and by part (i), $(F'_1,U'_1,U'_2,F'_2)$
  is a 3-walk. Suppose first that $\langle g'_1,g'_2\rangle$ acts
  irreducibly on $V$. Then by Lemma~\ref{L1}(b), we have $F'_1\cap F'_2= \{0\}$
  so $\{F'_1, F'_2\}$ is also an edge of $\Gamma_{e_1,e_2}$, so the
  3-walk $(F'_1,U'_1,U'_2,F'_2)$ is closed;
  also $U'_1\ne F'_2$ and $U'_2\ne F'_1$  by Lemma~\ref{L1}(c), and therefore
  $(F'_1,U'_1,U'_2,F'_2)$ is a 3-arc. Suppose conversely that
  $(F'_1,U'_1,U'_2,F'_2)$ is a closed 3-arc in $\Gamma_{e_1,e_2}$.  Then
  the vertices $F'_1, U'_1,U'_2, F'_2$ are pairwise distinct so
  Lemma~\ref{L1}(c) holds; $\{F'_1,F'_2\}$ is an edge so
  Lemma~\ref{L1}(b) holds; and also the condition $V=U'_1\oplus U'_2$ of
  Lemma~\ref{L1}(a) holds since $\{U'_1,U'_2\}$ is an edge. Thus
  $\langle g'_1,g'_2\rangle$ acts irreducibly by Lemma~\ref{L1}.
\end{proof}

We next prove that each of the maps $\phi_1,\phi_2, \psi$ in Lemma~\ref{L:irredclosed} has fibres of constant~size, and moreover that the $G_d$-conjugacy classes of stingray elements are also conjugacy classes for any group between $\SL_d(q)$ and $G_d$.

  For an $e$-stingray element $g\in\GL(V)$ let $U(g)\coloneq{\rm im}(g-1)$ and
  $F(g)\coloneq{\rm ker}(g-1)$.
  
\begin{lemma}\label{L:sting1}
  Let $g$ be an $e$-stingray element of $G_d=\GL_d(q)$. For any subgroup
  $G$ of $\GL_d(q)$ containing $\SL_d(q)$ let
  $\cC_g(G)=\{x^{-1}gx\mid x \in G\}$ be the $G$-conjugacy class of $g\in G_d$.
  We allow $g\not\in G$.
  \begin{enumerate}[{\rm (a)}]
  \item The $G$-conjugacy class $\cC_g(G)$ is independent of the choice of $G$, and has size
    \[
      |\cC_g(G)|=\frac{|G_d|}{(q^{e}-1)\cdot|G_{d-e}|}.
    \]
  \item The number of $g'\in\cC_g(G)$ such that $(U(g),F(g))=(U(g'),F(g'))$ is
    $|G_e|/(q^e-1)$.
  \item For $i=1,2$, let $g_i$ be an $e_i$-stingray element in $G_d$
    and let $\cC_i$ be the $G_d$-conjugacy class containing $g_i$.
    Then the number of pairs
    $(g'_1,g'_2)\in\cC_1\times\cC_2$ such that the $4$-tuple    $(F(g'_1),U(g'_1),U(g'_2),F(g'_2))$ equals\newline $(F(g_1),U(g_1),U(g_2),F(g_2))$, is 
  \begin{equation}\label{e:con}
  \frac{|G_{e_1}|\cdot |G_{e_2}|}{(q^{e_1}-1)(q^{e_2}-1)}.
  \end{equation}
  \end{enumerate}
\end{lemma}

\begin{proof}
  (a)~Let $h$ be the restriction of $g$
  to $U(g)$. Then $h$ is irreducible on $U(g)$ by
  Definition~\ref{d:sting}, and $C_{\GL(U(g))}(h)=Z_{q^{e}-1}$ is cyclic of order $q^{e}-1$
  by~\cite[Satz II.7.3, p.187]{Hup}.
  Therefore $X:=C_{G_d}(g)= Z_{q^{e}-1}\times \GL(F(G))$.
  In particular, $X$ contains a diagonal matrix of arbitrary non-zero
  determinant so $|\det(X)|=q-1$. Since $G$ contains $\SL_d(q)$, this implies that $G_d=GX$,~so
  $|\det(G)|$ divides $q-1$,~so
  \[
  |\cC_g(G)|=\frac{|G|}{|C_G(g)|}=\frac{|G|}{|G\cap X|}
  =\frac{|GX|}{|X|}=\frac{|G_d|}{|X|}=|\cC_g(G_d)|.
  \]
  As $\cC_g(G)\subseteq \cC_g(G_d)$, equality holds and $|\cC_g(G)|=\frac{|G_d|}{(q^{e}-1)|G_{d-e}|}$.
  
  (b)~By part (a), the number of $g'\in\cC_g(G)$ such that $U(g')=U(g)$ and $F(g')=F(g)$ is
  the number of $g'\in\cC_g(G_d)$ with this property. Thus we may,
  and shall, assume that $G=G_d$.
  Since $g$ is an $e$-stingray element, we have $V=U(g)\oplus F(g)$ with
  $\dim(U(g))=e$. If $g'=g^x$ where $x\in G$, we have $U(g')=U(g)x$ as
  \[
  U(g')=V(g'-1)= V(x^{-1}gx-1)=Vx^{-1}(g-1)x=V(g-1)x=U(g)x.
  \]
  Similarly $F(g')=F(g)x$.
  Thus the number of $g'\in\cC_g(G)$ such that
  $U(g')=U(g)$ and $F(g')=F(g)$ is equal to the number of choices for $x\in G$ with
  $U(g)x=U(g)$ and $F(g)x=F(g)$ divided by $|C_{G_d}(g)|$. As $G=G_d$, this
  equals $(|\GL(U(g))\times\GL(F(g))|)/|X|$ with $X=C_{G_d}(g)=Z_{q^{e}-1}\times \GL(F(G))$ as in part~(a). 
  Therefore, the number we seek is 
  $|\GL(U(g))|/(q^e-1) = |G_e|/(q^e-1)$, and part (b) is proved.
  
  (c) By part~(b) there are $|\GL_{e_1}(q)|/(q^{e_1}-1)$ choices for $g_1'\in\cC_1$ and $|\GL_{e_2}(q)|/(q^{e_2}-1)$ choices for $g_2'\in\cC_2$.
  Hence the number of pairs $(g'_1,g'_2)\in\cC_1\times\cC_2$ is 
  as claimed.
\end{proof}

\section{Closed 3-walks and closed 3-arcs of $\Gamma_{e_1,e_2}$}\label{S:3-walks}

In the light of Lemma~\ref{L1} we shall assume henceforth that
$e_1+e_2=d$ and therefore $V=\F_q^{e_1+e_2}$.  In this section we count the number of closed 3-walks and closed 3-arcs of~$\Gamma_{e_1,e_2}$. This is twice the number that occur in $\WW$ by Remark~\ref{R:WW}. Recall
from~\eqref{e:x12} that $X_i$ denotes the set of $e_i$-subspaces of $V=(\F_q)^{e_1+e_2}$ for $i\in\{1,2\}$.
It will be convenient to represent a subspace of $V$ as the row space of a block matrix.
Let $M_{e\times d}$ denote the vector space
of $e\times d$ matrices over $\F_q$. 
As the general linear group $G_d$ is transitive on the set of decompositions
$V=U_1\oplus U_2$ where $U_1\in X_1$ and $U_2\in X_2$, we write $U_1$ as the
row space of $\RSM{I}{0}\in M_{e_1\times d}$, and $U_2$ as the
row space of $\RSM{0}{I}\in M_{e_2\times d}$. We use the shorthand
$U_1=\RS{I}{0}$ and $U_2=\RS{0}{I}$ where the number of rows of $I$ and $0$
can be inferred from $\dim(U_i)=e_i$. The number of columns of $0$ can
be inferred, as $I$ is always a square matrix, and there are $d$ columns in
total. Therefore
$U_1=\RS{I}{0}$ has $0\in M_{e_1\times e_2}$ and $U_2=\RS{0}{I}$
has $0\in M_{e_2\times e_1}$.

We define the action of the group $G_{e_2}\times G_{e_1}$ on a matrix
  $A\in M_{e_2\times e_1}$, and on a pair
  $(A,B)\in M_{e_2\times e_1}\times M_{e_1\times e_2}$: for
  $(X,Y)\in G_{e_2}\times G_{e_1}$ write
  \begin{equation}\label{E:act}
    A^{(X,Y)}=X^{-1}AY\qquad\textup{and}\qquad
    (A,B)^{(X,Y)}=(X^{-1}AY,Y^{-1}BX).
  \end{equation}
  (These are actions as $(A^{(X_1,Y_1)})^{(X_2,Y_2)}=A^{(X_1X_2,Y_1Y_2)}$ and
  $((A,B)^{(X_1,Y_1)})^{(X_2,Y_2)}=(A,B)^{(X_1X_2,Y_1Y_2)}$.)

Since the
rank of $A$ equals the dimension of the row space of $A$, and
the dimension of the column space
of $A$, it follows that $A$ and $X^{-1}AY$
have the same rank. Indeed, the set of matrices of a given rank forms
an orbit under the group $G_{e_2}\times G_{e_1}$.
Since $e_2\leqslant e_1$, it follows that $G_{e_2}\times G_{e_1}$ has $e_2+1$
orbits on $M_{e_2\times e_1}$, and $G_{e_1}\times G_{e_2}$ has $e_2+1$
orbits on $M_{e_1\times e_2}$. Denote the set of $e_2\times e_1$ matrices
  rank~$k$ over $\F_q$ by $M_{e_2\times e_1}^{(k)}$.
  Let $\cR^{(k)}$ be 
$(G_{e_2}\times G_{e_1})$-orbit representatives for $M_{e_2\times e_1}^{(k)}$ where
\begin{equation}\label{E:Rk}
  \cR^{(k)}\coloneq\smatrix{I_k}{0}{0}{0}\in M_{e_2\times e_1}^{(k)}\subseteq M_{e_2\times e_1}
  \qquad\textup{for $0\leqslant k\leqslant e_2$.}
\end{equation}
The matrices in the second row of
$\cR^{(k)}$ have $e_2-k$ rows, and those in the 
second column have $e_1-k$ columns. In particular, $\cR^{(0)}$ is
the $e_2\times e_1$ zero matrix.

The cardinality $|M_{e_2\times e_1}^{(k)}|$ is
known, see Morrison~\cite[\S1.7]{Mo}. We give a short proof
to describe the structure of a stabiliser and to introduce our notation in the next lemma.

\begin{lemma}\label{L:OS}
  If $0\leqslant k\leqslant e_2\leqslant e_1$, then the stabiliser of $\cR^{(k)}$
  in $G_{e_2}\times G_{e_1}$ equals
  \[
  H_k=\left\{(X,Y)\in G_{e_2}\times G_{e_1}\mid
  X=\smatrix{X_{11}}{X_{12}}{0}{X_{22}}
  \textup{ and }
  Y=\smatrix{X_{11}}{0}{Y_{21}}{Y_{22}}\right\}\textup{ where}
  \]
  $X_{11}=Y_{11}\in G_k, X_{22}\in G_{e_2-k}, Y_{22}\in G_{e_1-k},
  X_{12}\in M_{k\times(e_2-k)}$ and $Y_{21}\in M_{(e_1-k)\times k}$.
  Hence the number of rank-$k$ matrices in $M_{e_2\times e_1}$ is
  $|M_{e_2\times e_1}^{(k)}|=|G_{e_1}\times G_{e_2}|/|H_k|$ where $|H_k|$ equals $|G_k||G_{e_1-k}||G_{e_2-k}|q^{k(e_1+e_2-2k)}$.
\end{lemma}

\begin{proof}
  The stabiliser $H_k$ of $\cR^{(k)}$ in $G_{e_2}\times G_{e_1}$ is easy to compute
  because the identity $\cR^{(k)}=(\cR^{(k)})^{(X,Y)}$ says that $X\cR^{(k)}=\cR^{(k)}Y$.
  In terms of matrices, this says
  \[
  \smatrix{X_{11}}{X_{12}}{X_{21}}{X_{22}} \smatrix{I_k}{0}{0}{0}
  =\smatrix{I_k}{0}{0}{0} \smatrix{Y_{11}}{Y_{12}}{Y_{21}}{Y_{22}}
  \quad\textup{that is}\quad
  \smatrix{X_{11}}{0}{X_{21}}{0}=\smatrix{Y_{11}}{Y_{12}}{0}{0}.
  \]
  Hence $X_{21}=0$, $Y_{12}=0$ and $X_{11}=Y_{11}\in G_k$.
  Thus $X_{22}\in G_{e_2-k}$ and
  $Y_{22}\in G_{e_1-k}$. Also $X_{12}\in M_{k\times(e_2-k)}$ and
  $Y_{21}\in M_{(e_1-k)\times k}$ are arbitrary. It follows that the
  stabiliser $H_k$ of $\cR^{(k)}$ is as claimed, and its order is
  $|G_k||G_{e_1-k}||G_{e_2-k}|q^{k(e_1+e_2-2k)}$. By the orbit-stabiliser lemma 
  the set of rank-$k$ matrices in $M_{e_2\times e_1}$ has
  cardinality  $|G_{e_2}\times G_{e_1}|/|H_k|$.
\end{proof}


When counting algebraic objects over $\F_q$, factoring out the
dominant power of $q$ determines the asymptotic behaviour as $q\to\infty$.
For example, $|\GL_e(q)|=q^{e^2}\omega(e)$ using \eqref{e:om}, and
$\omega(e)=\prod_{i=0}^{e-1}(1-q^{-i})=1-O(q^{-1})$
so $|\GL_e(q)|\sim q^{e^2}$ as $q\to\infty$.
By Lemma~\ref{L:OS}, $|M_{e_2\times e_1}^{(k)}|\sim q^x$
where $x=e_2^2+e_1^2-k^2-(e_1-k)^2-(e_2-k)^2-k(e_1+e_2-2k)=k(e_1+e_2-k)$.
A precise formula for $|M_{e_2\times e_1}^{(k)}|$ follows from Lemma~\ref{L:OS}:
\begin{equation}\label{E:Mk}
  |M_{e_2\times e_1}^{(k)}|=
  \frac{\omega(e_2)\omega(e_1)q^{k(e_1+e_2-k)}}{\omega(k)\omega(e_1-k)\omega(e_2-k)}
  \qquad\textup{where $0\leqslant k\leqslant e_2$}.
\end{equation}

\begin{corollary}\label{C}
  If $1\leqslant e_2\leqslant e_1$ and $q>1$ is a prime power, then using \eqref{e:om},
  \[
  \sum_{k=0}^{e_2}\frac{\omega(e_1)\omega(e_2)q^{-(e_1-k)(e_2-k)}}{\omega(k)\omega(e_1-k)\omega(e_2-k)}=
  \sum_{\ell=0}^{e_2}\frac{\omega(e_1)\omega(e_2)q^{-(e_1-e_2+\ell)\ell}}{\omega(e_2-\ell)\omega(e_1-e_2+\ell)\omega(\ell)}
  =1.
  \]
\end{corollary}

\begin{proof}
  Since $\sum_{k=0}^{e_2}|M_{e_2\times e_1}^{(k)}|=|M_{e_2\times e_1}|$, we have
  $\sum_{k=0}^{e_2}\frac{\omega(e_2)\omega(e_1)q^{k(e_1+e_2-k)}}{\omega(k)\omega(e_1-k)\omega(e_2-k)}=q^{e_1e_2}$ by~\eqref{E:Mk}. Setting $\ell=e_2-k$, the result
  follows from
  \[
  -e_1e_2+k(e_1+e_2-k)=-(e_1-k)(e_2-k)=-(e_1-e_2+\ell)\ell.
  \]
\end{proof}

We are now ready to count the closed $3$-walks.

\begin{theorem}\label{T:walk}
  Suppose that $d=e_1+e_2$ and $q>1$.
  Then the number of closed $3$-walks
  $(F_1,U_1,U_2,F_2)\in \WW$ in the bipartite
  $q$-Kneser graph $\Gamma_{e_1,e_2}$ equals
  \[
  q^{4e_1e_2}\omega(e_1+e_2)\sum_{\ell=0}^{e_2}
  \frac{q^{-(e_1-e_2+\ell)\ell}}
       {\omega(e_1-e_2+\ell)\omega(\ell)}.
  \]
\end{theorem}

\begin{proof}
As $G_d$ is transitive on decompositions $V=U_1\oplus U_2$, we may assume that
$U_1=\RS{I}{0}\in X_1$ and $U_2=\RS{0}{I}\in X_2$.
There are $q^{e_1e_2}$ complements $F_1\in X_2$ to $U_1$ in $V$ and these can
be written uniquely as $F_1=\RS{A}{I}$ where $A\in M_{e_2\times e_1}$ and
$I=I_{e_2}$ has $e_2$ rows (and columns). Similarly, the $q^{e_1e_2}$ complements
of $U_2$ can be written uniquely as $F_2=\RS{I}{B}\in X_1$ where
$B\in M_{e_1\times e_2}$. The number of pairs $(U_1,U_2)$ is
\[
  \frac{|G_{e_1+e_2}|}{|G_{e_1}\times G_{e_2}|}=q^{2e_1e_2}\xi \qquad
  \textup{where $\xi \coloneq\frac{\omega(e_1+e_2)}{\omega(e_1)\omega(e_2)}$.}
\]
Hence the number of
3-walks $(F_1,U_1,U_2,F_2)\in \WW$ is
$q^{4e_1e_2}\xi $ (as we saw in Lemma~\ref{L:number}). The number of closed
$3$-walks in $\WW$
is less than this, namely
$N\cdot q^{2e_1e_2}\xi $ where $N$ is the
number of choices of pairs $(F_1,F_2)$, given $(U_1,U_2)$, such that 
$V=U_1\oplus F_1=U_2\oplus F_2=F_1\oplus F_2$. We now find $N$.

The following are equivalent: $\{F_1,F_2\}$ is an edge of $\Gamma_{e_1,e_2}$;
$F_1\cap F_2$ equals $\{0\}$; and $F_1+ F_2$ equals $V$. In terms of matrices,
each of these conditions is equivalent to the constraint that the
$d\times d$ matrix $\smatrix{A}{I}{I}{B}$ is invertible. However,
$\smatrix{I}{-A}{0}{I}\smatrix{A}{I}{I}{B}=\smatrix{0}{I-AB}{I}{B}$. Therefore
the condition that $\{F_1,F_2\}$ lies in $E\Gamma_{e_1,e_2}$ is equivalent to
the condition that $I_{e_2}-AB$ is invertible. It is possible for
$I_{e_2}-AB$ to have full rank $e_2$, and lie in $G_{e_2}$, due to our
assumption that~$e_2\leqslant e_1$.

Thus we must determine the number $N$ of
pairs $(A,B)\in M_{e_2\times e_1}\times M_{e_1\times e_2}$ for which $I-AB$ is an invertible
$e_2\times e_2$ matrix. Recall the action~\eqref{E:act} of
$(X,Y)\in G_{e_2}\times G_{e_1}$ on a pair
$(A,B)$ in $M_{e_2\times e_1}\times M_{e_1\times e_2}$.
This action preserves pairs $(A,B)$ with $I-AB$ invertible since $(A',B')=(A,B)^{(X,Y)}$ implies
\[
  I-A'B'=I-(X^{-1}AY)(Y^{-1}BX)=I-X^{-1}ABX=X^{-1}(I-AB)X,
\]
and hence $\det(I-A'B')=\det(I-AB)$. Therefore $I-AB$ is invertible
precisely when $I-A'B'$ is invertible, and
\[
  N=\sum_{k=0}^{e_2}\#\{A\mid A\in M_{e_2\times e_1}^{(k)}\}
     \cdot \#\{B\in M_{e_1\times e_2}\mid I-AB\in G_{e_2}\}.
\]

For $A\in M_{e_2\times e_1}^{(k)}$, we may choose
$(X,Y)\in G_{e_2}\times G_{e_1}$ so that $X^{-1}AY=\cR^{(k)}$ is the rank-$k$
representative defined by~\eqref{E:Rk}. Therefore,
\[
  N=\sum_{k=0}^{e_2}|M_{e_2\times e_1}^{(k)}|
     \cdot \#\{B\in M_{e_1\times e_2}\mid I-\cR^{(k)}B\in G_{e_2}\}.
\]
Write $B=\smatrix{B_{11}}{B_{12}}{B_{21}}{B_{22}}$ where $B_{11}\in M_{k\times k}$.
Since $I-\cR^{(k)}B=\smatrix{I-B_{11}\,}{\,-B_{12}}{0}{I}$, a necessary and
sufficient condition for $I-\cR^{(k)}B$ to be invertible is that
$I-B_{11}\in G_k$. The matrices $B_{12},B_{21},B_{22}$ may be chosen
arbitrarily. However, the number of matrices $B_{11}$ with $I-B_{11}\in G_k$
is the number of $k\times k$ matrices not having 1 as an eigenvalue.
This equals the number of $k\times k$ matrices not having 0 as an eigenvalue,
that is $|G_k|=q^{k^2}\omega(k)$. In summary, the number of choices for
$B_{11},B_{12},B_{21},B_{22}$ is
$q^{k^2}\omega(k)$, $q^{k(e_2-k)}$, $q^{(e_1-k)k}$, $q^{(e_1-k)(e_2-k)}$, respectively.
Consequently, there are $q^{e_1e_2}\omega(k)$ choices for $B$.

The formula~\eqref{E:Mk} for $|M_{e_2\times e_1}^{(k)}|$ shows that
\[
N=\sum_{k=0}^{e_2}\frac{\omega(e_2)\omega(e_1)q^{k(e_1+e_2-k)}}{\omega(k)\omega(e_1-k)\omega(e_2-k)}\cdot q^{e_1e_2}\omega(k)
=\sum_{k=0}^{e_2}\frac{\omega(e_2)\omega(e_1)q^{e_1e_2+k(e_1+e_2-k)}}{\omega(e_1-k)\omega(e_2-k)}.
\]
Multiplying by the number $q^{2e_1e_2}\xi $ of pairs $(U_1,U_2)$, the number
of closed 3-walks is
\[
\frac{q^{2e_1e_2}\omega(e_1+e_2)}{\omega(e_1)\omega(e_2)}
\sum_{k=0}^{e_2}\frac{\omega(e_2)\omega(e_1)q^{e_1e_2+k(e_1+e_2-k)}}{\omega(e_1-k)\omega(e_2-k)}.
\]
Observing that $e_1e_2+k(e_1+e_2-k)=2e_1e_2-(e_1-k)(e_2-k)$, this equals
\[
  q^{4e_1e_2}\omega(e_1+e_2)\sum_{k=0}^{e_2}
  \frac{q^{-(e_1-k)(e_2-k)}}
       {\omega(e_1-k)\omega(e_2-k)}
  =q^{4e_1e_2}\omega(e_1+e_2)\sum_{\ell=0}^{e_2}
  \frac{q^{-(e_1-e_2+\ell)\ell}}
       {\omega(e_1-e_2+\ell)\omega(\ell)}
\]
where in the last step $\ell\coloneq e_2-k$ ranges from $e_2$ down to $0$.
\end{proof}

Now we determine the number of closed $3$-arcs.

\begin{theorem}\label{T:3arcs}
  Suppose that $1\leqslant e_2\leqslant e_1$ and $q>1$ is a prime power.
  Then the number of closed $3$-arcs
  $(F_1,U_1,U_2,F_2)\in \WW$ in the 
  $q$-Kneser graph $\Gamma_{e_1,e_2}$ equals
  \[
  q^{4e_1e_2}\omega(e_1+e_2)\sum_{\ell=0}^{e_2-1}
  \frac{q^{-(e_1-e_2+\ell)\ell}}{\omega(e_1-e_2+\ell)\omega(\ell)}
  \left(1-\frac{q^{-e_1e_2}}{\omega(e_2-\ell)}\right).
  \]
\end{theorem}

\begin{proof}
  Our proof uses the arguments and the notation of Theorem~\ref{T:walk}.
  We count the closed 3-arcs $(F_1,U_1,U_2,F_2)$
  with $V=U_1\oplus U_2$. There are
  $q^{2e_1e_2}\xi $ choices for the pair $(U_1,U_2)$ in $X_1\times X_2$. Each
  choice of $(U_1,U_2)$ gives rise to the same number
  of choices of
  $(F_1,F_2)$ in $X_2\times X_1$ such that $(F_1,U_1,U_2,F_2)$ is a closed 3-arc.
  Hence we take $U_1=\RS{I}{0}$ and $U_2=\RS{0}{I}$. Following the proof of
  Theorem~\ref{T:walk}, we write $F_1=\RS{A}{I}$ and $F_2=\RS{I}{B}$ where
  $A\in M_{e_2\times e_1}$ and $B\in M_{e_1\times e_2}$. We require that
  $F_1\ne U_2$, $F_2\ne U_1$, and that $I-AB\in M_{e_2\times e_2}$
  is invertible. That is, $A\ne0$, $B\ne0$ and $I-AB\in G_{e_2}$.

  Suppose that $A$ has rank $k$. As $A\ne0$ and $e_2\leqslant e_1$, we have
  $|M_{e_2\times e_1}^{(k)}|$ choices for $A$ where $k\in\{1,\dots,e_2\}$.
  Each $A\in M_{e_2\times e_1}^{(k)}$ gives the same number of choices for
  $B$, so we may assume that $A$ equals $R_{e_2\times e_1}^{(k)}$
  as per~\eqref{E:Rk}. Writing $B=\smatrix{B_{11}}{B_{12}}{B_{21}}{B_{22}}$
  with $B_{11}\in M_{k\times k}$, we require that $B\ne0$ and $I-B_{11}\in G_k$.
  If $B_{11}=0$, there are $q^{e_1e_2-k^2}-1$ choices for $B$, namely any
  non-zero matrix $B=\smatrix{0}{B_{12}}{B_{21}}{B_{22}}$. If $B_{11}\ne 0$, then
  $B_{11}=I-g$ where $g\in G_k$ is non-identity and $B_{12},B_{21},B_{22}$ are
  arbitrary, so the number of choices for $B$ is
  \[
    \left(q^{k^2}\omega(k)-1\right)
    \cdot q^{k(e_2-k)}\cdot q^{(e_1-k)k}\cdot q^{(e_1-k)(e_2-k)}
    =q^{e_1e_2}\left(\omega(k)-q^{-k^2}\right).
  \]
  Hence for each $A\in M_{e_2\times e_1}^{(k)}$ the total number of
  choices for $B\in M_{e_1\times e_2}$ is
  \[
  q^{e_1e_2}\omega(k)x_k\qquad\textup{where}\quad x_k\coloneq 1-\frac{q^{-e_1e_2}}{\omega(k)}.
  \]
  Therefore, the number of pairs $(F_1,F_2)$ is the number of pairs $(A,B)$ namely (using the formula \eqref{E:Mk})
  \[
  \sum_{k=1}^{e_2}\frac{\omega(e_2)\omega(e_1)q^{k(e_1+e_2-k)}}{\omega(k)\omega(e_1-k)\omega(e_2-k)}\cdot q^{e_1e_2}\omega(k)x_k
=\sum_{k=1}^{e_2}\frac{\omega(e_2)\omega(e_1)q^{e_1e_2+k(e_1+e_2-k)}x_k}{\omega(e_1-k)\omega(e_2-k)}.
\]
Multiplying by the number $q^{2e_1e_2}\xi $ of pairs $(U_1,U_2)$, the number of
closed 3-arcs is
\[
q^{2e_1e_2}\omega(e_1+e_2)
\sum_{k=1}^{e_2}\frac{q^{e_1e_2+k(e_1+e_2-k)}x_k}{\omega(e_1-k)\omega(e_2-k)}.
\]
Observing that $e_1e_2+k(e_1+e_2-k)=2e_1e_2-(e_1-k)(e_2-k)$, this equals
\[
  q^{4e_1e_2}\omega(e_1+e_2)\sum_{k=1}^{e_2}
    \frac{q^{-(e_1-k)(e_2-k)}x_k}{\omega(e_1-k)\omega(e_2-k)}
  =q^{4e_1e_2}\omega(e_1+e_2)\sum_{\ell=0}^{e_2-1}
  \frac{q^{-(e_1-e_2+\ell)\ell}x_{e_2-\ell}}
       {\omega(e_1-e_2+\ell)\omega(\ell)}
\]
where in the last step $\ell\coloneq e_2-k$ ranges from $e_2-1$ down to $0$.
\end{proof}

\subsection{Proof of Theorem \ref{TT:prop}: an explicit formula for $P(e_1,e_2)$}\label{sub:TTprop}

We finish this section by drawing together the results on $3$-walks and $3$-arcs in $\Gamma_{e_1,e_2}$ to prove Theorem~\ref{TT:prop}. Thus we assume that $q>1$ is a prime power and $d=e_1+e_2$ with
  $1\leqslant e_2\leqslant e_1$, and we consider the proportion $P(e_1,e_2)$  of $3$-walks of $\Gamma_{e_1,e_2}$ which are closed $3$-arcs. To prove Theorem~\ref{TT:prop} we must prove that \eqref{e:Pe12} holds for $P(e_1,e_2)$, that is, 
 \[
  P(e_1,e_2)=-(1-q^{-e_1e_2})q^{-e_1e_2}+\sum_{\ell=0}^{e_2-1}
  \frac{\omega(e_1)\omega(e_2)q^{-(e_1-e_2+\ell)\ell}}{\omega(e_1-e_2+\ell)\omega(\ell)}
  =1-O(1/q).
  \]
  In addition, it suffices to consider $3$-walks and $3$-arcs in $\WW$ by
    Remark~\ref{R:WW}.
    Let $V=(\F_q)^d$ and note that we may choose decompositions
    $V=(\F_q)^d=U_1\oplus F_1=U_2\oplus F_2$ with $U_1,F_2\in X_1$ of
    dimension $e_1$ and $U_2,F_1\in X_2$ of dimension $e_2$.
  
  The number $n$ of closed  3-arcs in $W_2\coloneq\WW$ is then given by Theorem~\ref{T:3arcs}, and the number of $3$-walks in $W_2$ is $q^{4e_1e_2}\xi $ by Lemma~\ref{L:number} where $\xi =\frac{\omega(e_1+e_2)}{\omega(e_1)\omega(e_2)}$. Thus by Remark~\ref{R:WW}, $P(e_1,e_2)=\frac{n}{q^{4e_1e_2}\xi }$ and we have
  \begin{align*}
  P(e_1,e_2)&=\frac{q^{4e_1e_2}\omega(e_1+e_2)}{q^{4e_1e_2}\xi }\sum_{\ell=0}^{e_2-1}
  \frac{q^{-(e_1-e_2+\ell)\ell}}{\omega(e_1-e_2+\ell)\omega(\ell)}
  \left(1-\frac{q^{-e_1e_2}}{\omega(e_2-\ell)}\right)\\
  &=\omega(e_1)\omega(e_2)\sum_{\ell=0}^{e_2-1}
  \frac{q^{-(e_1-e_2+\ell)\ell}}{\omega(e_1-e_2+\ell)\omega(\ell)}
  \left(1-\frac{q^{-e_1e_2}}{\omega(e_2-\ell)}\right).
  \end{align*}
  Corollary~\ref{C} implies that
  $\sum_{\ell=0}^{e_2-1}
  \frac{\omega(e_1)\omega(e_2)q^{-(e_1-e_2+\ell)\ell}}{\omega(e_2-\ell)\omega(e_1-e_2+\ell)\omega(\ell)}=1-q^{-e_1e_2}$. Hence
  \[
  P(e_1,e_2)=\left(\sum_{\ell=0}^{e_2-1}
  \frac{\omega(e_1)\omega(e_2)q^{-(e_1-e_2+\ell)\ell}}{\omega(e_1-e_2+\ell)\omega(\ell)}\right)
  -(1-q^{-e_1e_2})q^{-e_1e_2}
  \]
  as in \eqref{e:Pe12}. When $\ell=0$ the summand is $\frac{\omega(e_1)\omega(e_2)}{\omega(e_1-e_2)} = 1-O(1/q)$. Finally, this implies
  that $P(e_1,e_2)=1-O(1/q)$, as claimed.\qed

\section{Proportions: closed 3-arcs  and  stingray duos}\label{S:3-arcs}

The aim of this section is to prove Theorem~\ref{TT:stingray}. We do this in two steps. Our first result considers stingray duos in a fixed pair of conjugacy classes.

  \begin{theorem}\label{T:stingray}
  Let $d=e_1+e_2$ with  $1\leqslant e_2\leqslant e_1$,  
  let $\SL_d(q)\trianglelefteqslant G\leqslant \GL_d(q)$, and for $i=1,2$,    let $\cC_i$ be a $G$-conjugacy class of $e_i$-stingray elements. Then the proportion 
  \[
  \frac{\text{Number of irreducible stingray duos in  $\cC_1\times \cC_2$}}{\text{Number of stingray duos in  $\cC_1\times \cC_2$}}
  \]
  equals the proportion $P(e_1,e_2)$ given by \eqref{e:Pe12} of $3$-walks
  of $\Gamma_{e_1,e_2}$ which are closed $3$-arcs.
   \end{theorem}

 \begin{proof}
 By Lemma~\ref{L:sting1}(a), $\cC_1$  and $\cC_2$ are conjugacy classes of $G_d=\GL_d(q)$, and so we may assume that $G=G_d$. 
 It follows from  Lemma~\ref{L:sting1}(c) that, for each $3$-walk
 $(F_1,U_1,U_2,F_2)$ in $\WW$ of $\Gamma_{e_1,e_2}$, the equalities 
   \begin{equation}\label{e:ineq}
   \mbox{$U_i={\rm im}(g_i-1)$ and  $F_i={\rm ker}(g_i-1)$, for $i=1,2$,}
   \end{equation}
   all hold for the same number~\eqref{e:con} of stingray duos
   $(g_1,g_2)\in\cC_1\times\cC_2$.
   This is true, in particular, for closed $3$-arcs $(F_1,U_1,U_2,F_2)$. 
   Further, by Lemma~\ref{L:irredclosed}(b), either all or none of the
   pairs $(g_1,g_2)$ satisfying the equalities in \eqref{e:ineq} have
   $\langle g_1,g_2\rangle$ irreducible, and irreducibility occurs
   precisely when $(F_1,U_1,U_2,F_2)$ is a closed $3$-arc. Hence the
   proportion of $(e_1,e_2)$-stingray duos $(g_1,g_2)$ in
   $\cC_1\times\cC_2$ for which $\langle g_1,g_2\rangle$ is
   irreducible, equals the proportion of 3-walks
   $(F_1,U_1,U_2,F_2)\in\WW$ which are closed
   3-arcs.  Moreover, this proportion is the quantity $P(e_1, e_2)$ in
   \eqref{e:Pe12} by Remark~\ref{R:WW}.
 \end{proof}
 
 Theorem~\ref{T:stingray} is an important component in the proof of
 Theorem~\ref{TT:stingray} below.

\begin{proof}[Proof of Theorem~\ref{TT:stingray}]
Recall that $d=e_1+e_2$, $1\leqslant e_2\leqslant e_1$, and the
subgroup $G$ satisfies $\SL_d(q)\trianglelefteqslant G\leqslant
G_d=\GL_d(q)$.  Let $\cE_i$ be the set of all $e_i$-stingray elements
in $G$, for $i=1, 2$. Clearly $G$ acts on $\cE_i$. Let $\cY_i$ be
the set of $G$-conjugacy classes that partition~$\cE_i$. Then
$\cY_1\times\cY_2\coloneq\{\cC_1\times\cC_2\mid\cC_1\in\cY_1, \cC_2\in\cY_2\}$ is the set of $(G\times G)$-conjugacy classes
that partition $\cE_1\times\cE_2$. For a subset $Z$ of $\cE_1\times\cE_2$
let $D(Z)$ denote the set of $(e_1,e_2)$-stingray duos in $Z$,
and let $I(Z)$ denote the (sub)set of \emph{irreducible} $(e_1,e_2)$-stingray duos in $Z$.
Paraphrasing Theorem~\ref{T:stingray} gives
$|I(\cC_1\times\cC_2)|=P(e_1,e_2)\cdot |D(\cC_1\times\cC_2)|$. We must prove that
$|I(\cE_1\times\cE_2)|=P(e_1,e_2)\cdot |D(\cE_1\times\cE_2)|$. Since
$I(\cE_1\times\cE_2)$ is a disjoint union of $I(\cC_1\times\cC_2)$
where $\cC_1\times\cC_2$ ranges over $\cY_1\times\cY_2$, and similarly for
$D(\cE_1\times\cE_2)$ and $D(\cC_1\times\cC_2)$, we have
\begin{align*}
|I(\cE_1\times\cE_2)|
&=\sum_{\cC_1\times\cC_2}|I(\cC_1\times\cC_2)|
=\sum_{\cC_1\times\cC_2}P(e_1,e_2)\cdot |D(\cC_1\times\cC_2)|\\
&=P(e_1,e_2)\cdot |D(\cE_1\times\cE_2)|.
\end{align*}\qed
\end{proof}

\section{Explicit upper and lower bounds for  $P(e_1,e_2)$}\label{S:Order}

Recall that $P(e_1,e_2)$ is the proportion of $(e_1,e_2)$-stingray duos
in $\GL_d(q)$ that are irreducible by Theorem~\ref{TT:stingray}.
In this section we prove Theorem~\ref{TT:bounds} giving precise upper and
lower bounds for $P(e_1,e_2)$ where $1\leqslant e_2\leqslant e_1$, and
$q>1$ is an arbitrary prime power.

The ordering of two real numbers written in base-$q$ is determined by
the largest power of~$q$ where the digits differ. For example,
if $a(q)=\sum_{i\leqslant s} a_iq^i$ and $b(q)=\sum_{j\leqslant s} b_jq^j$ are (necessarily
convergent) Laurent series in $q$, with 
$a_i,b_j\in\{0,1,\dots,q-1\}$, $b_s\ne0$ and $a_s < b_s$,
then $a(q)\leqslant b(q)$ holds. Further, if $a_s<b_s$, then $a(q)=b(q)$ holds 
precisely when $b_s=a_s+1$ and,  for each $i< s$, $a_i=q-1$ and $b_i=0$.
An example with s = 0 is $\sum_{i<0} (q-1)q^i = q^0$.
Note that we can only compare Laurent series in $q$ if the coefficients are
not too large (in absolute value). For example, 
\[
1-2q^{-2}-2q^{-3}+q^{-4}+4q^{-5}+q^{-6}-2q^{-7}-2q^{-8}+q^{-10}\leqslant 1-2q^{-2}
\]
holds for all $q\geqslant5$. However, direct evaluation shows that it is also true
for $q\in\{2,3,4\}$.


\begin{proof}[Proof of Theorem~\ref{TT:bounds}]
Upper and lower bounds for $P(e_1,1)$ are best handled
separately. Substituting $e_2=1$ into the formula~\eqref{e:Pe12}
for $P(e_1,e_2)$ in Theorem~\ref{TT:prop} gives 
\begin{equation*}
    P(e_1,1)=-(1-q^{-e_1})q^{-e_1}+(1-q^{-1})(1-q^{-e_1})
   =(1-q^{-e_1})(1-q^{-1}-q^{-e_1}).
\end{equation*}
Therefore $P(e_1,1)=1-q^{-1}-2q^{-e_1}+q^{-e_1-1}+q^{-2e_1}$. For $e_1\geqslant3$
and for all~$q\geqslant2$ we have
$1-q^{-1}-q^{-2}< P(e_1,1)< 1-q^{-1}$.
(The upper bound even holds for $e_1\geqslant1$.)

Henceforth suppose that $2\leqslant e_2\leqslant e_1$.
The case when $e_1=e_2$ requires delicate estimations.
We first prove the upper bound when $e_1=e_2=e\geqslant2$.
  Theorem~\ref{TT:prop} gives
  \[
  P(e,e)= -(1-q^{-e^2})q^{-e^2}
  +\sum_{\ell=0}^{e-1}\frac{\omega(e)^2q^{-\ell^2}}{\omega(\ell)^2}.
  \]
  Setting $e=2$ gives
  \[
    P(2,2)=1-q^{-1}-q^{-2}+2q^{-3}-2q^{-4}-q^{-5}+q^{-6}+q^{-8}. 
  \]
  It follows that $1-q^{-1}-q^{-2}< P(2,2)< 1-q^{-1}-q^{-2}+2q^{-3}-2q^{-5}$.
  Suppose now  that $e_1=e_2=e\geqslant3$. Since $\omega(e)\leqslant\omega(3)$ and
  $\omega(\infty)<\omega(\ell)$, we have
  \begin{align*}
  P(e,e)&<\sum_{\ell=0}^{\infty}\frac{\omega(e)^2q^{-\ell^2}}{\omega(\ell)^2}
    <\omega(3)^2\left(1+\frac{q^{-1}}{\omega(1)^2}+\frac{q^{-4}}{\omega(2)^2}
  +\sum_{\ell=3}^{\infty}\frac{q^{-\ell^2}}{\omega(\infty)^2}\right)\\
  &<\omega(3)^2\left(1+\frac{q^{-1}}{\omega(1)^2}+\frac{q^{-4}}{\omega(2)^2}
  +\sum_{\ell=3}^{\infty}\frac{q^{-8(\ell-2)}}{\omega(\infty)^2}\right);
  \end{align*}
  where the last line uses the fact that $8(\ell-2)\leqslant \ell^2$ for all $\ell\geqslant3$.
  Adding the geometric series $\sum_{\ell=3}^{\infty}q^{-8(\ell-2)}$ gives
  $\frac{q^{-8}}{1-q^{-8}}$. We use the inequality
  $1-q^{-1}-q^{-2}+q^{-5}<\omega(\infty)$, which follows from the proof
  of~\cite[Lemma~3.5]{NP},  so that
  \begin{equation}\label{E:xyz}
  P(e,e)<\omega(3)^2\left(1+\frac{q^{-1}}{\omega(1)^2}+\frac{q^{-4}}{\omega(2)^2}
  +\frac{q^{-8}}{(1-q^{-1}-q^{-2}+q^{-5})^2(1-q^{-8})}\right).
  \end{equation}
  Since $(1-q^{-1}-q^{-2}+q^{-5})^2(1-q^{-8})> \frac{1}{16}\geqslant q^{-4}$
  for all $q$, we have by~\eqref{E:xyz} that
  \begin{align*}
  P(e,e)&<\omega(3)^2\left((1+q^{-4})+\frac{q^{-1}}{\omega(1)^2}
    +\frac{q^{-4}}{\omega(2)^2}\right)\\
  &=\omega(3)^2(1+q^{-4})+(1-q^{-2})^2(1-q^{-3})^2q^{-1}+(1-q^{-3})^2q^{-4}.
  \end{align*}
  The inequality below is true for $q\geqslant5$ by comparing Laurent series
  in $q$:
  \begin{align*}
  (1-q^{-2})^2(1-q^{-3})^2
    &=1-2q^{-2}-2q^{-3}+q^{-4}+4q^{-5}+q^{-6}-2q^{-7}-2q^{-8}+q^{-10}\\
    &\leqslant 1-2q^{-2}.
  \end{align*}
  However, the inequality $(1-q^{-2})^2(1-q^{-3})^2\leqslant 1-2q^{-2}$ is also true for
  $q\in\{2,3,4\}$. Hence $\omega(3)^2<(1-q^{-1})^2(1-2q^{-2})$, and so
  \begin{align*}
  P(e,e)&<(1-q^{-1})^2(1-2q^{-2})(1+q^{-4})+(1-2q^{-2})q^{-1}+(1-q^{-3})^2q^{-4}\\
  &=1-q^{-1}-q^{-2}+2q^{-3}-2q^{-5}-q^{-6}+2q^{-7}-2q^{-8}+q^{-10}.
  \end{align*}
  The above reasoning proves that the
  upper bound $P(e,e)<1-q^{-1}-q^{-2}+2q^{-3}-2q^{-5}$
  holds for all $e\geqslant2$ and all prime powers $q\geqslant2$.

  We next prove that the lower bound $1-q^{-1}-q^{-2}< P(e,e)$ holds for
  $e\geqslant3$ and $q\geqslant2$. Setting $e_1=e_2=e$ in
  Theorem~\ref{TT:prop} gives
  \begin{align*}
    P(e,e)&= -(1-q^{-e^2})q^{-e^2}+\sum_{\ell=0}^{e-1}\frac{\omega(e)^2q^{-\ell^2}}{\omega(\ell)^2}
    > -q^{-9}+\omega(e)^2\sum_{\ell=0}^{2}\frac{q^{-\ell^2}}{\omega(\ell)^2}.
  \end{align*}
  Now
  \begin{align*}
  \omega(2)^2\sum_{\ell=0}^{2}\frac{q^{-\ell^2}}{\omega(\ell)^2}
  &=1-q^{-1}-q^{-2}+2q^{-3}-q^{-5}+q^{-6}\\
  &>1-q^{-1}-q^{-2}+2q^{-3}-q^{-5}=(1-q^{-2})(1-q^{-1}+q^{-3}).\\
  \end{align*}
  Using  $\omega(e)>\omega(\infty)>1-q^{-1}-q^{-2}+q^{-5}$ and the above inequality gives
  \begin{align*}
    P(e,e)&>-q^{-9}+\frac{\omega(e)^2}{\omega(2)^2}(1-q^{-2})(1-q^{-1}+q^{-3})
    =-q^{-9}+\frac{\omega(e)^2(1-q^{-1}+q^{-3})}{(1-q^{-1})^2(1-q^{-2})}\\
    &>-q^{-9}+\frac{(1-q^{-1}-q^{-2}+q^{-5})^2(1-q^{-1}+q^{-3})}{(1-q^{-1})^2(1-q^{-2})}\kern-0.8pt.
  \end{align*}
  We approximate the denominator using:
  $(1-q^{-1})^{-2}=\sum_{i=0}^\infty(i+1)q^{-i}>\sum_{i=0}^6(i+1)q^{-i}$ and
  $(1-q^{-2})^{-1}=\sum_{i=0}^\infty q^{-2i}>\sum_{i=0}^2 q^{-2i}$. Hence
  \def\pm{\kern-1pt-\kern-1pt}
  \def\pp{\kern-1pt+\kern-1pt}
   \begin{align*}
    P(&e,e)\kern-1pt>\kern-1pt -q^{-9}+(1-q^{-1}-q^{-2}+q^{-5})^2(1-q^{-1}+q^{-3})
    \left(\sum_{i=0}^6(i+1)q^{-i}\right)\left(\sum_{i=0}^2q^{-2i}\right)\\
    &=1\pm q^{-1}\pm q^{-2}\pp q^{-5}\pm2q^{-6}\pm8q^{-7}\pp16q^{-8}\pp q^{-9}\pm5q^{-10}\pp q^{-11}\pm5q^{-12}\pp26q^{-13}\\
    &\quad \pm12q^{-14}\pm7q^{-15}\pm q^{-16}\pm9q^{-17}\pp6q^{-18}\pm17q^{-19}\pm3q^{-20}\pp5q^{-21}\pp6q^{-22}\pp7q^{-23}.
  \end{align*}
  The above Laurent series is therefore greater than $1-q^{-1}-q^{-2}$
  for $q\geqslant27$, and a computer calculation shows it is greater
  than $1-q^{-1}-q^{-2}$ for $3\leqslant q<26$. The same lower bound for $P(e,e)$ also holds for $q=2$ if we
  replace the bound $\omega(\infty)>1-q^{-1}-q^{-2}+q^{-5}=0.28125$ with
  the sharper inequality $\omega(\infty)>0.288$.
  Thus the bound $P(e,e)>1-q^{-1}-q^{-2}$ holds
  for $e\geqslant3$ and all prime powers $q\geqslant2$.
  This establishes the lower bound in Theorem~\ref{TT:bounds} for the case $e_1=e_2$.
  Henceforth assume that $e_1\ne e_2$.

  We next prove that $P(e_1,e_2)< 1-q^{-1}-q^{-2}+2q^{-3}-2q^{-5}$ when $2\leqslant e_2<e_1$.
  It is convenient to write
  $P(e_1,e_2)=-(1-q^{-e_1e_2})q^{-e_1e_2}+L(e_1,e_2)+M(e_1,e_2)$ where
  \[
  L(e_1,e_2)\coloneq\sum_{\ell=0}^{1}
  \frac{\omega(e_1)\omega(e_2)q^{-(e_1-e_2+\ell)\ell}}{\omega(e_1-e_2+\ell)\omega(\ell)}\hskip2mm \textup{and}\hskip2mm
  M(e_1,e_2)\coloneq\sum_{\ell=2}^{e_2-1}
  \frac{\omega(e_1)\omega(e_2)q^{-(e_1-e_2+\ell)\ell}}{\omega(e_1-e_2+\ell)\omega(\ell)}.
  \]
  When $e_2$ or $q$ is large, the term $L(e_1,e_2)$ dominates $P(e_1,e_2)$.
  Rearranging gives
\begin{align}
  L(e_1,e_2)&\coloneq\sum_{\ell=0}^{1}
  \frac{\omega(e_1)\omega(e_2)q^{-(e_1-e_2+\ell)\ell}}{\omega(e_1-e_2+\ell)\omega(\ell)}=
  \frac{\omega(e_1)\omega(e_2)}{\omega(e_1-e_2)}
+\frac{\omega(e_1)\omega(e_2)q^{-(e_1-e_2+1)}}{\omega(e_1-e_2+1)\omega(1)}\notag\\
&=\frac{\omega(e_1)\omega(e_2)}{\omega(e_1-e_2+1)\omega(1)}
  \left((1-q^{-(e_1-e_2+1)})(1-q^{-1})+q^{-(e_1-e_2+1)}\right)\notag\\
&=\frac{\omega(e_1)\omega(e_2)}{\omega(e_1-e_2+1)\omega(1)}\label{EEE}
  \left(1-q^{-1}+q^{-(e_1-e_2+2)}\right).
\end{align}
The inequalities $\omega(e_1)<\omega(e_1-e_2+\ell)$,
$\omega(e_2)<\omega(\ell)$ and $e_1-e_2\geqslant1$ imply that
\[
M(e_1,e_2)<\sum_{\ell=2}^{e_2-1}  q^{-(\ell+1)\ell}
<q^{-6}+\sum_{\ell=3}^\infty  q^{-(\ell+1)\ell}<q^{-6}+\sum_{\ell=3}^{\infty}  q^{-4\ell}\leqslant q^{-6}+2q^{-12}\kern-0.9pt<\kern-0.9pt q^{-6}+q^{-7}.
\]

If $e_2=2$, then $M(e_1,e_2)=0$. Using~\eqref{EEE} and the inequality
$\omega(e_1)<\omega(e_1-e_2+1)$~gives
\begin{align*}
P(e_1,2)&<L(e_1,2) =(1-q^{-e_1})(1-q^{-2})(1-q^{-1}+q^{-e_1})\\
&=(1-q^{-2e_1}-q^{-1}+q^{-(e_1+1)})(1-q^{-2})\\
&<(1-q^{-1}+q^{-4})(1-q^{-2})=1-q^{-1}-q^{-2}+q^{-3}+q^{-4}-q^{-6}\\
&<1-q^{-1}-q^{-2}+2q^{-3}-2q^{-5}.
\end{align*}
Suppose now that $3\leqslant e_2<e_1$. In this case
$\omega(e_2)\leqslant\omega(3)$, so we have the upper bounds
$L(e_1,e_2)<\frac{\omega(3)}{\omega(1)}(1-q^{-1}+q^{-3})$
and $M(e_1,e_2)<q^{-6}+q^{-7}$. Hence
\begin{align*}
  P(e_1,e_2)&<L(e_1,e_2)+M(e_1,e_2)
  < \frac{\omega(3)}{\omega(1)}(1-q^{-1}+q^{-3})+q^{-6}+q^{-7}\\
&=(1-q^{-2})(1-q^{-3})(1-q^{-1}+q^{-3})+q^{-6}+q^{-7}\\
  &=1-q^{-1}-q^{-2}+q^{-3}+q^{-4}-q^{-6}+q^{-7}+q^{-8}\\
  &< 1-q^{-1}-q^{-2}+2q^{-3}-2q^{-5}.
\end{align*}
Thus the bound $P(e_1,e_2)< 1-q^{-1}-q^{-2}+2q^{-3}-2q^{-5}$ holds for 
$2\leqslant e_2\leqslant e_1$ and $q\geqslant2$.

We now prove that $1-q^{-1}-q^{-2}< P(e_1,e_2)$ for
$2\leqslant e_2< e_1$.
By Theorem~\ref{TT:prop},
$P(e_1,e_2)=-(1-q^{-e_1e_2})q^{-e_1e_2}+L(e_1,e_2)+M(e_1,e_2)$.
If $e_2=2$, then $M(e_1,2)=0$ and
\begin{align*}
  P(e_1,2)&=-(1-q^{-2e_1})q^{-2e_1}+L(e_1,2)\\
  &=-(1-q^{-2e_1})q^{-2e_1}+(1-q^{-e_1})(1-q^{-2})(1-q^{-1}+q^{-e_1})
  \quad\quad{\rm by~\eqref{EEE}}\\
  &> -q^{-2e_1}+(1-q^{-2})(1-q^{-1}+q^{-e_1}-q^{-e_1}+q^{-e_1-1}-q^{-2e_1})\\
  &>-q^{-4}+(1-q^{-2})(1-q^{-1})\\
  &=-q^{-4}+1-q^{-1}-q^{-2}+q^{-3} > 1-q^{-1}-q^{-2}.
\end{align*}
Suppose now that $3\leqslant e_2< e_1$.
Since $-(1-q^{-e_1e_2})q^{-e_1e_2}\geqslant-(1-q^{-9})q^{-9}> -q^{-9}$ and
$M(e_1,e_2)\geqslant0$, we have $P(e_1,e_2)> -q^{-9}+L(e_1,e_2)$.
We next find a sharp
lower bound for $L(e_1,e_2)$ for $3\leqslant e_2< e_1$.
Set $\ell_0\coloneq e_1-e_2+2$. Then $3\leqslant\ell_0< e_1$ and
\[
\frac{\omega(e_1)}{\omega(\ell_0-1)}=\prod_{i=\ell_0}^{e_1}(1-q^{-i})
> 1-\sum_{i=\ell_0}^{e_1} q^{-i} >1-\sum_{i=\ell_0}^\infty q^{-i}
=1-\frac{q^{-\ell_0}}{1-q^{-1}}.
\]
In addition,
\[
  \frac{\omega(e_2)}{\omega(1)} >\frac{\omega(\infty)}{1-q^{-1}}.
\]
Multiplying the previous two inequalities shows:
\begin{equation}\label{E:twice}
  \frac{\omega(\infty)(1-q^{-1}-q^{-\ell_0})}{(1-q^{-1})^2}
  < \frac{\omega(e_1)\omega(e_2)}{\omega(e_1-e_2+1)\omega(1)}.
\end{equation}
Recall that $\ell_0\coloneq e_1-e_2+2$ and $3\leqslant e_2< e_1$. Using the
formula~\eqref{EEE} for $L(e_1,e_2)$ and equation~\eqref{E:twice} gives
\begin{align*}
  P(e_1,e_2)&>-q^{-9}
    +\frac{\omega(\infty)(1-q^{-1}-q^{-\ell_0})(1-q^{-1}+q^{-\ell_0})}{(1-q^{-1})^2}\\
  &=-q^{-9}+\frac{\omega(\infty)(1-2q^{-1}+q^{-2}-q^{-2\ell_0})}{(1-q^{-1})^2}.
\end{align*}
We next use $\omega(\infty)>1-q^{-1}-q^{-2}+q^{-5}$ (see \cite[Lemma~3.5]{NP}) and
\[
  \frac{1}{(1-q^{-1})^2}=\sum_{i=0}^\infty(i+1)q^{-i}> \sum_{i=0}^6(i+1)q^{-i}.
\]
Since $\ell_0\geqslant3$, we have
\begin{align*}
  P(e_1,e_2)&> -q^{-9} +(1-q^{-1}-q^{-2}+q^{-5})(1-2q^{-1}+q^{-2}-q^{-6})
    \left(\sum_{i=0}^6(i+1)q^{-i}\right)\\
  &=1-q^{-1}-q^{-2}+q^{-5}-q^{-6}-9q^{-7}+15q^{-8}+q^{-9}-5q^{-10}\\
  &\qquad+2q^{-11}-6q^{-12}+17q^{-13}+3q^{-14}-5q^{-15}-6q^{-16}-7q^{-17}\\
  &\geqslant 1-q^{-1}-q^{-2};
\end{align*}
the last inequality holds since
$q^{-5}-q^{-6}-9q^{-7}+15q^{-8}+\cdots-7q^{-17}>0$ for $q>17$,
and direct calculation shows that it also holds for $2\leqslant q\leqslant17$.
This establishes the lower bound of Theorem~\ref{TT:bounds} for $2\leqslant e_2< e_1$ and $q\geqslant2$, and completes
the proof.
\end{proof}

\section{Motivation from Computational Group Theory}\label{S:comp}

We noted in Section~\ref{sec:intro} that the results of this paper relate to
a larger problem of recognising finite classical groups,
see~\cite{DLLOB,GNP2022a,GNP2022b,LOB,PSY}. This section
describes the context of this problem, and why the lower bound $1-q^{-1}-q^{-2}$
for the proportion $P(e_1,e_2)$
given in this paper is helpful, while the lower bound $1-2q^{-1}+O(q^{-2})$
given in~\cite{PSY}, for a related proportion, is problematic.

Let $\GX_d(q)$ denote a classical group over the field $\F_{q^\delta}^d$
of type $\mathbf{X}$, where $\mathbf{X}$ is $\mathbf{L}$ (linear), $\mathbf{U}$ (unitary), $\mathbf{S}$ (symplectic), or $\mathbf{O}$ (orthogonal) and
$\delta=1$ unless $\mathbf{X}=\mathbf{U}$ when $\delta=2$.
Then $\GX_d(q)$ acts naturally on a vector space $V=\F_{q^\delta}^d$.
If $\GX_d(q)$ is not solvable, let $\OmegaX_d(q)$ be the smallest normal
subgroup of $\GX_d(q)$ for which $\GX_d(q)/\OmegaX_d(q)$ is solvable
and let $\SX_d(q)$ denote the subgroup of $\GX_d(q)$
comprising  matrices with determinant~$1$.
If $\mathbf{X}=\mathbf{L},\mathbf{S},\mathbf{U}$
then $\SX_d(q)=\OmegaX_d(q)$, whereas $|\SO_d(q):\Omega_d(q)|=2$ when
$\mathbf{X}=\mathbf{O}$.
We present the key strategy for constructive recognition algorithms of classical groups, first introduced in the
state-of-the-art algorithm \cite{DLLOB}, a strategy which is also adopted in
a new algorithm currently being developed.

Conceptually we start with a group $G=\langle A\rangle$ which is
(either known to be, or believed to be) a classical group $G=\SY_n(q)$  of type $\mathbf{Y}$ acting on a vector space $V=\F_{q^\delta}^n$ as above.
The aim of a constructive recognition algorithm is to write the elements of a pre-defined generating set for $\SY_n(q)$ as words in $A$.
Algorithms following the strategy 
in \cite{DLLOB}
are recursive. In the first step, 
see for example \cite[p.\ 232 or Section 5]{DLLOB}, they
 construct a `first subgroup $H$' which,
with respect to an appropriate basis, has the form 
\[
H=   \smatrix{\SX_{d}(q)}{0}{0}{I_{n-d}} 
\]
where  $d < n$ and either  $\mathbf{X}=\mathbf{Y}$ or  $(\mathbf{Y},\mathbf{X},q)=(\mathbf{S},\mathbf{O},\mbox{even})$.
To construct $H$, the algorithms seek a 
\emph{generating}
$(e_1,e_2)$-stingray duo, namely a
$(e_1,e_2)$-stingray duo
$(g_1,g_2)$ in $\GX_d(q)$ for which 
 $\langle g_1,g_2\rangle$ contains $\OmegaX_d(q)$.
The
correctness and the complexity analysis, as well as the practical runtime,
of a given constructive recognition algorithm, is underpinned by
tight estimates for the proportion 
of $(e_1,e_2)$-stingray duos $(g_1,g_2)$
in $\OmegaX_d(q)$ (or in some specified pair of $\OmegaX_d(q)$-conjugacy classes) that are generating and for which each $g_i$ has prime order.
For example, the algorithm in  \cite{DLLOB}
seeks pairs of $(e_1,e_2)$-stingray duos of
the form $(g,g^x)$ with $g\in G$
and $g^x$ a random $G$-conjugate,
where 
$e=e_1=e_2=d/2$ for some $d\in  [n/3, 2n/3].$
With high probability $(g, g^x)$ is a generating
 $(e,e)$-stingray duo, see \cite[Section 5.1]{DLLOB}.
This is implied by  \cite[Lemma 5.8]{DLLOB} in case $\mathbf{L}$, and in  \cite[Lemma 3.4]{GNP2022b} for cases  $\textbf{X}=\textbf{S},\textbf{U},$ and $\textbf{O}$. (There was an unfortunate gap in extending the proof of \cite[Lemma 5.8]{DLLOB} to cases other than $\textbf{L}$ which was repaired in \cite[Lemma 3.4]{GNP2022b}, see \cite[Remark 3.5]{GNP2022b}.)
In forthcoming algorithms, we choose $d$ to be $O(\log(n))$
and construct $(e_1,e_2)$-stingray duos with $d=e_1+e_2.$

The results in the present paper yield substantial improvements for
the analysis of constructive recognition algorithms in three ways.

Firstly, they are required for the analysis
of the above mentioned forthcoming algorithm, which hinges on proving that 
the proportion of generating $(e_1,e_2)$-stingray duos among 
all $(e_1,e_2)$-stingray duos in $\OmegaX_{e_1+e_2}(q)$ is large, which will be
proved in~\cite{GNP2023}.
 
Secondly, they improve existing estimates for the proportion 
of $(e,e)$-stingray duos of the form $(g,g^x)$ for $g\in
\OmegaX_d(q)$ that are generating,
where  difficulties in handling small values of
$q$ have resulted  in restrictions in the  applicability of complexity
estimates for these algorithms. For  example, the main theorem for the
constructive recognition  of classical  groups in  even characteristic
\cite[Theorem 1.2]{DLLOB} requires $q>4$ because the estimation result
about  stingray elements  it relied  on  (from a  preprint version  of
\cite{PSY})   had  this   restriction.   In   the  published   version
\cite[Theorem  2]{PSY},   the  only  restriction  is   `$q>2$  in  the
orthogonal    case',    see    our    discussion    in    \cite[Remark
  3.5(b)]{GNP2022b}.

Finally, our results improve the bounds for the complexity analysis of existing algorithms,
by noting that the proportions used to estimate the performance of algorithms
underestimated the proportion of elements actually constructed in the algorithm.

We address these aspects in the following two subsections.

\subsection*{Estimating  the proportion  of $(e_1,e_2)$-stingray duos}

Let  $G$  be   a  group  satisfying  $\OmegaX_d(q)\trianglelefteqslant
G\leqslant\GX_d(q)$.   The  proportion  of  $(e_1,e_2)$-stingray  duos
$(g_1,g_2)$  in  $G$  (or  in some  specified  pair  of  $G$-conjugacy
classes)  that are  generating  can be  estimated  by bounding  the
 proportion  of  $(e_1,e_2)$-stingray  duos $(g_1,g_2)$ which do not generate, and therefore lie in a proper 
 maximal    subgroup   of    $\OmegaX_{d}(q).$
This approach was pioneered in \cite[Theorem 6]{PSY}
for a  classical group $\SX_{d}(q)$  of type $\mathbf{X}$  with $d=2e$
even,  and  an  $\SX_{d}(q)$-conjugacy  class  $\cC$  of  $e$-stingray
elements of prime order:  the results  \cite[Theorem  2, 5,  6]{PSY}  show that  with
probability $1-O(q^{-\delta})$ a pair in $\cC\times \cC$ will generate
a  subgroup $\SX'_{d}(q)$  of  type $\mathbf{X}'$  (provided
$q>2$ if $\mathbf{X}=\mathbf{O}$), where as noted above,  $\mathbf{X}'=\mathbf{X}$ or  $(\mathbf{X},\mathbf{X}',q)=(\mathbf{S},\mathbf{O},\mbox{even})$.
Equivalently, the probability that
such  a pair  $(g,g^x)$ fails  to generate  is $O(q^{-\delta})$.   The
probability that $\langle  g,g^x\rangle$ is irreducible and  lies in a
proper maximal subgroup  of $\SX_{d}(q)$, i.e.\ it  fails to generate,
is shown  to be  very small, namely  $O(q^{-cd^2})$ where  $c$ is  a constant
depending on the type $\mathbf{X}$, \cite[Theorem 6]{PSY}. Therefore,
the leading
terms in  the estimates  correspond to pairs  where $\langle
g,g^x\rangle$  is reducible  on  $\F_{q^{\delta}}^d$.  In  forthcoming
work \cite{GNP2023} we are finding  that a similar dichotomy holds for
general  $(e_1,  e_2)$-stingray duos.  Also  here,   the  most
difficult problem arises  when $\mathbf{X}=\mathbf{L}$, where
the probability that an $(e_1,e_2)$-stingray duo is not generating is dominated by  the
probability that an $(e_1,e_2)$-stingray duo is reducible.
 Using the lower bound for $P(e_1,e_2)$  obtained in Theorem~\ref{TT:bounds}  we arrive
at such a bound.

\begin{corollary}\label{C:red}
  Suppose that $d=e_1+e_2$ with $2\leqslant e_2\leqslant e_1$, and
  $\SL_d(q)\lhdeq G\leqslant\GL_d(q)$. Then
  \[
  \frac{\text{Number of reducible $(e_1,e_2)$-stingray duos in $G\times G$}}
       {\text{Number of $(e_1,e_2)$-stingray duos in $G\times G$ }} \le q^{-1} + q^{-2}.
\]
Similarly, if  $\cC_i$ is a  $\GL_d(q)$-conjugacy class of $e_i$-stingray
elements for $i=1,2$, then
    \[
  \frac{\text{Number of reducible $(e_1,e_2)$-stingray duos in $\cC_1\times \cC_2$}}
       {\text{Number of $(e_1,e_2)$-stingray duos in $\cC_1\times \cC_2$ }} \le q^{-1} + q^{-2}.
\]
\end{corollary}

\begin{proof}
  The stated proportion is $1-P(e_1,e_2)$ by Theorem~\ref{T:stingray}. Then by Theorem~\ref{TT:bounds},
  $1-q^{-1}-q^{-2}<P(e_1,e_2)$ holds, so that $1-P(e_1,e_2)<q^{-1}+q^{-2}$ as claimed. The bound when restricted to stingray duos in $\cC_1\times \cC_2$ follows by the same argument applying Theorem~\ref{TT:stingray} rather than Theorem~\ref{T:stingray}.
\end{proof}

\subsection*{Stingray pairs versus stingray duos}

We now discuss how our results improve the analysis 
of the algorithm given in \cite{DLLOB} originally based on the results in \cite{PSY}. 
The bounds did not accurately reflect the practical performance 
in the linear case $\mathbf{L}$ where the classical group  is $\SL_d(q)$ with natural module $\F_q^d$, where $d=2e$ is even.  Let $\cC$ be an  $\SL_{d}(q)$-conjugacy  class  of
$e$-stingray elements.
By \cite[Lemma 5.3]{PSY},  the proportion of pairs  from $\cC\times \cC$ that generate a reducible subgroup is at most
\begin{equation}\label{e:q2}
    2q^{-1}+q^{-2}-2q^{-3}-q^{-4}+2q^{-d^2/4}.
\end{equation}
By  \cite[Theorem~6]{PSY}\footnote{In the course of our work for \cite{GNP2023} we discovered a problem with the proof of \cite[Lemma 10.1]{PSY} which is corrected and generalised in \cite{GNP2023}, so \cite[Theorem 6]{PSY} is valid.},
the proportion of pairs from $\cC\times \cC$ that generate an irreducible
proper subgroup of $\SL_d(q)$ is 
$O(q^{-d^2/4+d/2+2})$.
Thus, for sufficiently large $d$, the probability that a random pair  $\cC\times \cC$  generates $\SL_d(q)$ is at least $1-2q^{-1}+O(q^{-2})$, which 
is not a very  useful bound if $q=2$!

Our suggested solution is to exploit the fact that the algorithm in \cite{DLLOB}, and other algorithms under development, do not choose random pairs in $\cC\times \cC$, rather they construct stingray duos. Thus the appropriate proportions to estimate are of the form:
\[
\frac{\text{Number of stingray duos in $\cC\times \cC$ with a desired property}}{\text{Number of stingray duos in $\cC\times \cC$ }},
\]
or in the more general setting where we are considering $(e_1,e_2)$-stingray duos coming from $\cC_1\times \cC_2$, where $\cC_i$ is a conjugacy class of $e_i$-stingray elements, estimates of the form:
\[
\frac{\text{Number of stingray duos in $\cC_1\times \cC_2$ with a desired property}}{\text{Number of stingray duos in $\cC_1\times \cC_2$ }}.
\]
In other words, the proportions required for the complexity analysis of the recognition algorithms mentioned above 
require a \emph{different  denominator}. 
For the case where the desired property we wish to estimate is reducibility, Corollary~\ref{C:red} yields
an upper bound that is independent of $e_1, e_2$ and, in particular, is an
improvement on the bound obtained by applying Equation~\eqref{e:q2} to
the case relevant to \cite{DLLOB, PSY}.
However, to derive estimates for the proportion of stingray duos which
are irreducible but non-generating using the estimates
from~\cite[Theorem 6]{PSY}, we need to know the proportion of pairs
in $\cC_1\times \cC_2$ which are stingray duos. 
This proportion is given in Theorem~\ref{T:duo}~below.

\begin{theorem}\label{T:duo}
  Let $d=e_1+e_2$  and  let $\cC_i$ be a $\GL_d(q)$-conjugacy class of
  $e_i$-stingray elements, for $i=1,2$. Then 
 \[
\frac{\text{Number of stingray duos in $\cC_1\times \cC_2$ }}{|\cC_1\times\cC_2|} 
= \frac{1}{\xi},  \qquad\textup{where $\xi=\frac{\omega(d)}{\omega(e_1)\omega(e_2)}$}.
\]
Moreover, if $2\leqslant e_2\leqslant e_1$, then 
 $\frac{(1-q^{-d})(1-q^{-(d-1)})}{(1-q^{-1})(1-q^{-2})}\leqslant \xi<\frac{1}{1-q^{-1}-q^{-2}+q^{-5}}$.
\end{theorem}
    
\begin{proof} Let $G_d=\GL_d(q)$. 
    By Lemma~\ref{L:sting1}(a), 
    \[
    |\cC_1\times\cC_2| = \frac{|G_d|}{|G_{e_2}|\cdot(q^{e_1}-1)}\cdot  \frac{|G_d|}{|G_{e_1}|\cdot(q^{e_2}-1)} = \frac{|G_d|^2}{|G_{e_1}|\cdot|G_{e_2}|\cdot (q^{e_1}-1)(q^{e_2}-1)} 
    \]
    and we note that 
    \[
    \frac{|G_d|}{|G_{e_1}|\cdot|G_{e_2}|} = q^{d^2-e_1^2-e_2^2} \xi = q^{2e_2e_2}\xi.
    \]
    Next, by Lemmas~\ref{L:irredclosed}(b)(ii) and~\ref{L:sting1}(c), the number of stingray duos in  $\cC_1\times \cC_2$ is equal to the number of $3$-walks in $\Gamma_{e_1,e_2}$ times $\frac{|G_{e_1}|\cdot|G_{e_2}|}{(q^{e_1}-1)(q^{e_2}-1)}$. Moreover,  by Lemma~\ref{L:number}, the number of $3$-walks in  $\Gamma_{e_1,e_2}$ is $q^{4e_1e_2}\xi$. Thus the number $N$ of stingrays duos in  $\cC_1\times \cC_2$ is equal to 
    \[
 N=   q^{4e_1e_2}\xi \cdot \frac{|G_{e_1}|\cdot|G_{e_2}|}{(q^{e_1}-1)(q^{e_2}-1)}.
    \]
    The  proportion of stingray duos in
  $\cC_1\times \cC_2$ is therefore equal to 
    \[
  \frac{N}{|\cC_1\times \cC_2|} =    q^{4e_1e_2}\xi \cdot \frac{|G_{e_1}|\cdot|G_{e_2}|}{(q^{e_1}-1)(q^{e_2}-1)} \cdot 
    \frac{|G_{e_1}|\cdot|G_{e_2}|\cdot (q^{e_1}-1)(q^{e_2}-1)}{|G_d|^2} =\frac{1}{\xi}
    \]
    as claimed. In~\cite{GNP2022b} the $q$-binomial notation
   $\binom{d}{e}_q=\frac{\omega(d)}{\omega(e)\omega(d-e)}$ is used. Using this notation $\xi=\binom{d}{e_2}_q$, and 
   $\frac{(1-q^{-d})(1-q^{-(d-1)})}{(1-q^{-1})(1-q^{-2})}\leqslant \xi< \omega(\infty)^{-1}$   by~\cite[Lemma~5.2]{GNP2022b}, with $\omega(\infty)$ as in \eqref{e:om}. Moreover  $\omega(\infty)> 1-q^{-1}-q^{-2} + q^{-5}$   by~\cite[Lemma~5.1]{GNP2022b}.
\end{proof}

Note that Theorem~\ref{T:duo} yields a new upper bound on the proportion of pairs from $\cC_1\times \cC_2$ that generate a reducible subgroup. 

\begin{lemma}\label{C:red2}
  Let $d=e_1+e_2$ with $2\leqslant e_2\leqslant e_1$,  and  let $\cC_i$ be a $\GL_d(q)$-conjugacy class of
  $e_i$-stingray elements, for $i=1,2$. Then the proportion of pairs from $\cC_1\times \cC_2$ that generate a reducible subgroup is $1-\frac{P(e_1,e_2)}{\xi}$ with $\xi$ as in Theorem~\ref{T:duo}, and 
   \[
  1-\frac{P(e_1,e_2)}{\xi} < 2q^{-1} + q^{-2} - 2q^{-3} - q^{-4}.
 \]   
\end{lemma}

\begin{proof}
    Let $\mathcal{I}$ be the set of pairs $(g_1,g_2)\in\cC_1\times \cC_2$ such that $\langle g_1,g_2\rangle$ is irreducible, so the proportion we need to bound is $1-|\mathcal{I}|/|\cC_1\times \cC_2|$. If   $(g_1,g_2)\in  \cC_1\times \cC_2$ and $\langle g_1,g_2\rangle$ is irreducible, then the corresponding subspaces $U_1, U_2$ in Definition~\ref{d:duo} must be disjoint, as otherwise $\langle g_1,g_2\rangle$ would leave invariant the proper subspace $U_1+U_2$. Hence $(g_1,g_2)$ is a stingray duo. Thus $\mathcal{I}$ is the set of irreducible stingray duos in $\cC_1\times \cC_2$. Therefore, if $N= \text{Number of stingray duos in $\cC_1\times \cC_2$}$, then Theorem~\ref{T:stingray} gives
    \[
   \frac{|\mathcal{I}|}{|\cC_1\times \cC_2|} = \frac{|\mathcal{I}|}{N}\cdot \frac{N}{|\cC_1\times\cC_2|} = P(e_1,e_2)\cdot \frac{N}{|\cC_1\times\cC_2|}. 
    \]
    Theorem~\ref{T:duo} gives $\frac{N}{|\cC_1\times\cC_2|}=\frac{1}{\xi}> 1-q^{-1}-q^{-2}+q^{-5}$ and, in addition, Theorem~\ref{TT:bounds} gives $P(e_1,e_2)> 1- q^{-1} - q^{-2}$. The following calculation completes the proof as
    $q^{-6} + q^{-7}<q^{-5}$
 \begin{align*}
 1- \frac{|\mathcal{I}|}{|\cC_1\times \cC_2|}&= 1 - \frac{P(e_1,e_2)}{\xi} < 1-  (1- q^{-1} - q^{-2})\cdot (1-q^{-1}-q^{-2}+q^{-5}) \\
 &=2q^{-1} +q^{-2} - 2q^{-3} - q^{-4} - q^{-5} + q^{-6} + q^{-7}.
 \end{align*}  
\end{proof}

The upper bound in Lemma~\ref{C:red2} with $e_1=e_2$ always improves the bound
in \eqref{e:q2}. When $q=2$, this yields the upper bound $15/16$ for all $d$.
Finally, the following lower bound for the proportion of pairs sought in
the algorithm in \cite{DLLOB} improves its analysis, and better 
reflects the fact that the algorithm produces stingray duos.

\begin{corollary}\label{C:DLLOB}
Let $\cC$ be an  $\SX_{d}(q)$-conjugacy  class  of
$e$-stingray elements with $e=\frac{d}{2}$. Then, for some positive constant $c$,
\[
\frac{\textrm{Number of generating stingray duos in $\cC\times \cC$}}{\textrm{Number of stingray duos in $\cC\times \cC$ }} > 1- q^{-1}-q^{-2}- c\cdot q^{-d^2/4+d/2+2}.
\]
\end{corollary}
\begin{proof}
  Applying the bound in Theorem~\ref{T:duo} together with \cite[Theorem 6]{PSY}
  to  the linear case $\mathbf{L}$ where $\cC_1=\cC_2=\cC$ is an $\SL_d(q)$-conjugacy class of $e$-stingray elements for $e=d/2$ of order equal to a
  primitive prime divisor of $q^{d/2}-1$ we find:
\begin{align*}
   &{} \frac{\textrm{Number of irreducible but non-generating stingray duos in $\cC\times \cC$}}{\textrm{Number of stingray duos in $\cC\times \cC$ }}=\\
    &{}\frac{\textrm{Number of irreducible but non-generating stingray duos in $\cC\times \cC$}}{|\cC\times \cC|/\xi} < c\cdot q^{-d^2/4+d/2+2} 
\end{align*}
for some positive constant $c$. By Theorem~\ref{TT:stingray}, the proportion of irreducible $(e,e)$-stingray duos is $P(e,e)$, and by Theorem~\ref{TT:bounds}, $P(e,e)> 1- q^{-1}-q^{-2}$, and hence  the proportion of stingray duos in $\cC\times \cC$ which  generate $\SL_d(q)$ is at least
\[
P(e,e)- c \cdot q^{-d^2/4+d/2+2} > 1- q^{-1}-q^{-2}- c\cdot q^{-d^2/4+d/2+2}.
\]
\end{proof}



\section{Acknowledgements}
The authors gratefully acknowledge support from the ARC Research Council Discovery Project Grant DP190100450. The second author
  acknowledges funding by the Deutsche Forschungsgemeinschaft (DFG, German Research Foundation) - Project 286237555 - TRR 195. Finally, we thank the referee for suggesting improvements to the paper.



\begin{thebibliography}{99}

\bibitem{Biggs}
  Norman Biggs,
  Algebraic graph theory.
  Second edition. Cambridge Mathematical Library. Cambridge University Press,
  Cambridge, 1993. viii+205 pp. 

\bibitem{Brouwer2010}
  A.\,E. Brouwer,
  The eigenvalues of oppositeness graphs in buildings of spherical type,
  In {\it Combinatorics and graphs}, volume {\bf 531} of Contemp. Math., 
  pp. 1--10. Amer. Math. Soc., Providence, RI, 2010.

\bibitem{DMM}
  Jan De Beule, Sam Mattheus and Klaus Metsch,
  An algebraic approach to Erd\H{o}s-Ko-Rado sets of flags in
  spherical buildings.
  \emph{J. Combin. Theory Ser. A} {\bf192} (2022), paper no. 105657, 33 pp.


\bibitem{DLLOB}
  Heiko Dietrich, C.\,R. Leedham-Green, Frank L\"{u}beck and E.\,A. O'Brien,
  Constructive recognition of classical groups in even characteristic.
  \emph{J. Algebra} {\bf391} (2013), 227--255.

\bibitem{GIM}
  S.P. Glasby, Ferdinand Ihringer and Sam Mattheus.
  The proportion of non-degenerate complementary subspaces in classical spaces. 
  \emph{Des. Codes, Cryptogr.} {\bf 91(9)} (2023) 2879--2891.  
  
\bibitem{GNP2022a}
  S.\,P. Glasby, Alice C. Niemeyer and Cheryl E. Praeger,
  The probability of spanning a classical space by two
  non-degenerate subspaces of complementary dimensions,
  \emph{Finite Fields Their Appl.} {\bf 82} (2022), paper no. 102055.

\bibitem{GNP2022b}
  S.\,P. Glasby, Alice C. Niemeyer and Cheryl E. Praeger,
  Random generation of direct sums of finite non-degenerate subspaces,
  \emph{Linear Algebra Appl.} {\bf 649} (2022), 408--432.

\bibitem{GNP2023}
  S.\,P. Glasby, Alice C. Niemeyer and Cheryl E. Praeger,
  The probability that two elements with large 1-eigenspaces generate
  a classical group (in prep.).

\bibitem{Hup}
  B. Huppert, \emph{Endliche Gruppen I}. Springer-Verlag, Berlin, 1967.


\bibitem{KS}
  William M. Kantor and \'{A}kos Seress, 
  Black box classical groups,
  \emph{Mem. Amer. Math. Soc.} {\bf149} (2001), no. 708, viii+168 pp. 

\bibitem{LOB}
  C.\,R. Leedham-Green and E.\,A. O'Brien,
  Constructive recognition of classical groups in odd characteristic.
  \emph{J. Algebra} {\bf322} (2009), no. 3, 833--881.

\bibitem{Mo}
  Kent E. Morrison, Integer sequences and matrices over finite fields.
  \emph{J. Integer Seq.} {\bf 9} (2006), no. 2, Article 06.2.1, 28 pp.

\bibitem{NP}
  Peter M. Neumann and Cheryl E. Praeger,
  Cyclic matrices over finite fields.
  \emph{J. London Math. Soc.} (2) {\bf52} (1995), no. 2, 263--284. 

\bibitem{PSY}
  Cheryl E. Praeger, \'{A}kos Seress, \c{S}\"{u}kr\"{u} Yal\c{c}{\i}nkaya,
  Generation of finite classical groups by pairs of elements with
  large fixed point spaces. \emph{J. Algebra} {\bf 421} (2015), 56--101.

\bibitem{DR}
Daniel Rademacher, 
Constructive recognition of finite classical groups with stingray elements. \emph{Dissertation, RWTH Aachen University} (2024).
\url{https://publications.rwth-aachen.de/record/995067}
\end{thebibliography}
\end{document}